\documentclass{article}

\pdfoutput=1

\usepackage{arxiv}
\usepackage[utf8]{inputenc} % allow utf-8 input
\usepackage[T1]{fontenc}    % use 8-bit T1 fonts
\usepackage{hyperref}       % hyperlinks
\usepackage{url}            % simple URL typesetting
\usepackage{booktabs}       % professional-quality tables
\usepackage{amsfonts}       % blackboard math symbols
\usepackage{nicefrac}       % compact symbols for 1/2, etc.
\usepackage{microtype}      % microtypography
\usepackage{graphicx}
\usepackage{natbib}
\usepackage{multirow}
\usepackage{makecell}
\usepackage[normalem]{ulem}
\usepackage{doi}
\usepackage{amsmath,amssymb,amsfonts}
\usepackage{algorithm,algorithmic}
\usepackage{amsthm}

\newtheorem{theorem}{\indent Theorem}[section]
\newtheorem{lemma}[theorem]{\indent Lemma}

\newtheorem{definition}{\indent Definition}[section]

\newtheorem{remark}{\indent Remark}[section] 

\title{A note on the renormalization of floating-point expansions}

\author{ 
	Yuanyuan Yang \\
	yuanyuanyang@whu.edu.cn \\
}

\renewcommand{\shorttitle}{A note on the renormalization of floating-point expansions}

\hypersetup{
	pdftitle={A note on the renormalization of floating-point expansions},
	pdfsubject={cs.MS, math.NA},
	pdfauthor={Yuanyuan Yang},
	pdfkeywords={multiple precision, floating-point arithmetic, Fast2Sum, floating-point expansions},
}

\begin{document}
	\maketitle
	
\begin{abstract} 
Based on the newly established condition for the validity of the Fast2Sum algorithm, we present a detailed analysis of the behavior of the VecSum algorithm for various input configurations. Furthermore, we prove the correctness of the quad-double addition algorithm implemented in the QD library and demonstrate that, in many practical cases, 2Sum can be safely replaced by Fast2Sum.
\end{abstract}
	
% keywords can be removed

\section{Introduction}
Double-double (DD), triple-double (TD), and quad-double (QD) arithmetic represent a real number as the unevaluated sum of two, three, or four binary64 floating-point numbers conforming to the IEEE 754 standard \cite{IEEE7542008}. TD and QD are special cases of floating-point expansions (FPEs) \cite{JM2016CAMPARY,Priest1991FPExpansion,RichardShewchuk1997}, in which a real number is represented as the unevaluated sum of $n$ floating-point numbers.

DD, TD, and QD are useful when a precision higher than binary64 is required, for example, in numerically sensitive semidefinite programming \cite{Nakata2010SDPQDDD}, in the study of chaotic dynamical systems such as finding sinks of the H\'{e}non map \cite{HenonMapGPU2014}, and in experimental mathematics \cite{Bailey2015HPMathPhy,Bailey2012HPinMathPhyandDyn}. Implementations of multiple-precision arithmetic fall into two categories: multi-term and multi-digit. The QD library of Bailey et al. \cite{Bailey2000QD} is representative of the former, and MPFR \cite{MPFR2007} of the latter. Multi-digit libraries such as MPFR have the advantage that the precision and the rounding behavior can be arbitrary, whereas multi-term libraries can usually guarantee only error bounds. Hence a multi-digit library is preferable when precise control over precision and rounding behavior is required. When performance is the main concern, multi-term arithmetic is preferred, since it can exploit the hardware floating-point unit.

However, renormalization of an FPE is much more complicated than that of a double-double. Two basic building blocks of FPE renormalization are the VecSum and VecSumErrBranch algorithms. Here we use the new condition for the validity of the Fast2Sum algorithm (Lemma \ref{lem:fast2sumextra}) to review the VecSum algorithm (Algorithm \ref{algo:vecsum}) and Algorithm \ref{algo:QDadd} of the QD library \cite{Bailey2000QD}.

\section{Notation and Preliminaries}
In this paper, we work in a radix-$2$, precision-$p$ floating-point (FP) arithmetic system, which we assume to have an unlimited exponent range, thereby precluding any occurrence of underflow or overflow. Here, $p$ denotes the precision of the FP arithmetic system, and $\mathbb{F}_p$ denotes the set of all floating-point numbers of precision $p$.

Our analysis focuses primarily on the round-to-nearest (RN) rounding mode. Throughout this paper, $\mathrm{RN}(x)$ denotes the floating-point number closest to $x$, with ties broken by the ties-to-even rule. The unit roundoff, denoted by $\mathbf{u}$, is defined as $\mathbf{u} = 2^{-p}$ for RN rounding.

For $x \in \mathbb{F}_p$, the notations $\mathrm{ufp}(x)$, $\mathrm{ulp}(x)$, and $\mathrm{uls}(x)$ denote the unit in the first place, unit in the last place, and unit in the last significant place, respectively. The symbol $e_x$ denotes the exponent of the FP number $x$. For instance, for a double-precision FP number $x\in\mathbb{F}_{53}$ written as
\begin{equation}
x=(1.x_1x_2\cdots x_{49}100)_2\cdot 2^e = (1x_1x_2\cdots x_{49}100)_2 \cdot 2^{e-52},  
\end{equation}
we have $\mathrm{ufp}(x)=2^{e}$, $\mathrm{ulp}(x)=2^{e-52}$, $\mathrm{uls}(x)=2^{e-50}$, and $e_x=e$.

It is straightforward to verify that  
\begin{equation*}
2^{-p}|x|\leq\mathrm{ulp}(x) \leq 2^{-p+1}|x|, |x-\mathrm{RN}(x)|\leq \mathbf{u}|x|.
\end{equation*}

Two fundamental error-free transformations, Fast2Sum and 2Sum, are described below. The Fast2Sum algorithm was originally introduced by Dekker \cite{Dekker1971Split}.
\begin{algorithm}[H]
	\algsetup{linenosize=\small}
	\small
	\caption{Fast2Sum}\label{algo:Fast2Sum}
	\textbf{Input:} $a, b$. The floating-point exponents $e_a$ and $e_b$ satisfy $e_a\ge  e_b$\\
	\textbf{Output:} $ (s,t)$ s.t. $s+t=a+b$ 
	\begin{algorithmic}
		\STATE $ s \leftarrow \mathrm{RN}(a+ b)$
		\STATE $ t \leftarrow \mathrm{RN}( b-\mathrm{RN}(s-a) )$ 
	\end{algorithmic}
\end{algorithm}

The conventional 2Sum algorithm, due to Knuth \cite{TAOCP} and M\o{}ller \cite{10.1007/BF01975722}, does not require $e_a\geq e_b$.	
\begin{algorithm}[H]
	\algsetup{linenosize=\small}
	\small
	\caption{2Sum}\label{algo:2Sum}
	\textbf{Input:} $ a,b $\\
	\textbf{Output:} $ (s,t)$ s.t. $s+t=a+b$
	\begin{algorithmic}
		\STATE $    s \leftarrow  \mathrm{RN}(a+b)$
		\STATE $	a' \leftarrow \mathrm{RN}(s-b) $
		\STATE $	b' \leftarrow \mathrm{RN}(s-a') $
		\STATE $ 	\delta_a \leftarrow \mathrm{RN}(a-a') $
		\STATE $ 	\delta_b \leftarrow \mathrm{RN}(b-b')$
		\STATE $	t \leftarrow 	\mathrm{RN}(\delta_a+\delta_b)$
	\end{algorithmic}
\end{algorithm}

\begin{lemma}\label{lem:fast2sumextra}
	Assume that the floating-point numbers $a, b \in \mathbb{F}_p$ satisfy $\mathrm{uls}(a) \geq \mathrm{ulp}(b)$. Then $s + t = a + b$ holds exactly. 
\end{lemma}
In previous work \cite{yang2024robustnessdoublewordadditionalgorithms}, we established that the essential condition for the validity of Fast2Sum is $\mathrm{uls}(a)\geq \mathrm{ulp}(b)$ under faithful rounding, a result later rediscovered in \cite{park2026odderrorfreefasttwosumgeneral}. This lemma can also be extended to radix-3 floating-point numbers. In \cite{yang2024robustnessdoublewordadditionalgorithms} we exploited this property to analyze in detail the behavior of double-word addition algorithms; here we apply the same property to study the VecSum algorithm and the QDAdd algorithm (Algorithm \ref{algo:QDadd}).

An FPE represents a real number $x$ as the unevaluated sum of several floating-point numbers: $x=\sum_{i=0}^{n-1}x_i$. To achieve higher precision and facilitate manipulation, we require that the significant bits, or at least the nonzero significant bits, of adjacent terms $x_i$ and $x_{i+1}$  be separated by a sufficient gap. In this case, we say that the FPE $\{x_i\}_{i=0}^{n-1}$ is nonoverlapping.

Several definitions of nonoverlapping have been proposed in the literature, for example, $\mathcal{B}$-nonoverlapping by Bailey \cite{930115}, $\mathcal{P}$-nonoverlapping by Priest \cite{PriestOnProperties}, ulp-nonoverlapping in \cite{Popescu2017Phd}, F-nonoverlapping by Fabiano et al. \cite{JMuller2019TripleDouble}, and $\mathcal{S}$-nonoverlapping by Shewchuk \cite{RichardShewchuk1997}. The formal definitions of these conditions are as follows.
\begin{definition}[$\mathcal{B}$-nonoverlapping]
	An FPE $x=\{x_i\}_{i=0}^{n-1}$ is $\mathcal{B}$-nonoverlapping if for all $0<i<n$, $|x_i|\leq \frac 1 2\mathrm{ulp}(x_{i-1})$.
\end{definition}
\begin{definition}[$\mathcal{P}$-nonoverlapping]
	An FPE $x=\{x_i\}_{i=0}^{n-1}$ is $\mathcal{P}$-nonoverlapping if for all $0<i<n$, $|x_i|< \mathrm{ulp}(x_{i-1})$.
\end{definition}
\begin{definition}[ulp-nonoverlapping]
	An FPE $x=\{x_i\}_{i=0}^{n-1}$ is ulp-nonoverlapping if for all $0<i<n$, $|x_i|\leq \mathrm{ulp}(x_{i-1})$.
\end{definition}
\begin{definition}[F-nonoverlapping]
	An FPE $x=\{x_i\}_{i=0}^{n-1}$ is F-nonoverlapping if for all $0<i<n$, $|x_i|\leq \frac 1 2 \mathrm{uls}(x_{i-1})$.
\end{definition}	
\begin{definition}[$\mathcal{S}$-nonoverlapping]
	An FPE $x=\{x_i\}_{i=0}^{n-1}$ is $\mathcal{S}$-nonoverlapping if for all $0<i<n$, $|x_i|< \mathrm{uls}(x_{i-1})$.
\end{definition}

An expansion may contain interleaving zeros, but the definitions above apply only to the non-zero terms of the expansion.

Usually, the intermediate result of an FPE arithmetic algorithm is an unevaluated sum $\{y_i\}_{i=0}^{m-1}$ with $n\leq m$. However, the expansion $\{y_i\}_{i=0}^{m-1}$ is no longer nonoverlapping, and only much weaker assumptions can be made about it. Hence a renormalization step is necessary. 

VecSum was first introduced by Priest \cite{Priest1991FPExpansion}. The name ``VecSum'' was coined by Ogita et al. \cite{Rump2005Kfold}. The VecSum algorithm, presented in Algorithm \ref{algo:vecsum}, is a chain of 2Sum operations that performs error-free transformations on the inputs, ensuring that $\sum_{i=0}^{m-1}y_i=\sum_{i=0}^{m-1}x_i$. If the inputs satisfy any of the conditions below, $\{y_i\}_{i=0}^{m-1}$ is an $\mathcal{S}$-nonoverlapping expansion with possible interleaving zeros \cite{Popescu2017Phd}. 
\begin{enumerate}
	\item $\{x_i\}_{i=0}^{m-1}$ overlaps by at most $d\leq p-2$ bits, with possible interleaving zeros, and $\frac {2^d} {1-2^{d-p}} (1+(m-2)2^{-p})\leq 2^{p-1}$; the latter always holds in practice. Moreover, the 2Sum in VecSum can be replaced by Fast2Sum for $p\geq 6$ and $m \leq 34$.  
	\item $\{x_i\}_{i=0}^{m-1}$ overlaps by at most $d\leq p-1$ bits, with possible interleaving zeros, provided that $m\leq 12$ for binary32 and $m\leq 39$ for binary64.
	\item $\{x_i\}_{i=0}^{m-2}$ is $\mathcal{S}$-nonoverlapping with possible interleaving zeros, and $x_{m-1}$ is a random floating-point number.
\end{enumerate}

Another interesting condition is that $\{x_i\}_{i=0}^{m-1}$ overlaps by at most $d\leq p-2$ bits with possible interleaving zeros and may contain pairs of at most two consecutive terms that overlap in $p$ digits. Under this condition, the output $\{y_i\}_{i=0}^{m-1}$ satisfies $|y_0|> |y_1|\geq \cdots > |y_{i-1}|\geq |y_i|>|y_{i+1}\geq |y_{i+2}|>\cdots$ and  $ |y_j|<\mathrm{uls}(y_{j-1})$ or  $|y_j|,|y_{j-1}|$ are equal to a power of $2$ \cite{Popescu2017Phd}. 

We now give an informal definition of ``overlap by at most $d$ digits'' as used in the above conditions. That is, if the terms in $\{x_i\}_{i=0}^{m-1}$ are written in positional notation, the digits of any two successive nonzero terms coincide in at most $d$ digit positions, and no three terms coincide in the same digit position \cite{Priest1991FPExpansion}.

\begin{algorithm}[H]
	\algsetup{linenosize=\small}
	\small
	\caption{VecSum$(x_0,x_1,\cdots,x_{m-1})$}\label{algo:vecsum}
	\begin{algorithmic}[1]
		\STATE $t_{m-1} \leftarrow x_{m-1}$
		\FOR {$i\leftarrow m-2$ \TO $0$} 
		\STATE $(t_i,y_{i+1})\leftarrow \mathrm{2Sum}(x_i,t_{i+1})$
		\ENDFOR
		\STATE $y_0\leftarrow t_0$
		\RETURN  $y_0,y_1,\cdots,y_{m-1}$
	\end{algorithmic}
\end{algorithm}

The $\mathcal{S}$-nonoverlapping condition is relatively weak. For instance, the significant bits of a floating-point number with precision $2np$ may be distributed across $3n$ components $\{x_i\}_{i=0}^{3n-1}$, each of precision $p$. Thus, an $\mathcal{S}$-nonoverlapping expansion of length $n$ does not necessarily provide a total precision of $np$. But if we gather and condense the significant bits in an $\mathcal{S}$-nonoverlapping array and remove the interleaving zeros, we may shorten the array without altering the sum. This is exactly how VecSumErrBranch (Algorithm \ref{algo:vseb}) works. 

\begin{algorithm}[H]
	\algsetup{linenosize=\small}
	\small
	\caption{VecSumErrBranch$(y_0,y_1,\cdots,y_{m-1},n)$\cite{930115,Priest1991FPExpansion}}\label{algo:vseb}
	\begin{algorithmic}[1]
		\STATE $j\leftarrow 0$
		\STATE $ \varepsilon_0\leftarrow y_0$
		\FOR{ $ i \leftarrow 0 $ \TO $m-2$}
		\STATE $ (z_j,\varepsilon_{i+1})\leftarrow \mathrm{Fast2Sum}(\varepsilon_{i},y_{i+1})$
		\IF {$ \varepsilon_{i+1}\neq 0 $}
		\IF {$j \geq n-1$}
		\RETURN $z_0,z_1,\cdots,z_{n-1}$ 
		\ENDIF
		\STATE $j\leftarrow j+1$
		\ELSE 
		\STATE $\varepsilon_{i+1}\leftarrow z_j$
		\ENDIF
		\ENDFOR
		\IF {$\varepsilon_{m-1} \neq 0$ \AND $j<n$ }
		\STATE $z_j \leftarrow \varepsilon_{m-1}$
		\STATE $j\leftarrow j+1$ 
		\ENDIF
		\STATE $z_{j,\cdots, n-1} \leftarrow 0$
		\RETURN  $z_0,z_1,\cdots,z_{n-1}$   
	\end{algorithmic}
\end{algorithm}

VecSumErrBranch resembles VecSum in structure but processes the input terms in descending order of magnitude. If the inputs of VecSumErrBranch come from VecSum, then, provided that no underflow or overflow occurs during the computation, $\{z_i\}_{i=0}^{n-1}$ is ulp-nonoverlapping. Since $(y_0,y_1)$ is the result of the last 2Sum in VecSum, the first iteration $i=0$ can be skipped \cite{Popescu2017Phd}.

If $(x_0,x_1,\cdots,x_{m-1})$ satisfies one of the conditions for VecSum, the combination of VecSum and VecSumErrBranch (Algorithm \ref{algo:renrm}) renders an overlapping array ulp-nonoverlapping and truncates it to $n$ terms.

\begin{algorithm}[H]
	\algsetup{linenosize=\small}
	\small\caption{Renormalize$(x_0,x_1,\cdots,x_{m-1},n)$}\label{algo:renrm}
	\begin{algorithmic}[1]
		\STATE $y[0:m-1]\leftarrow \mathrm{VecSum}(x[0:m-1])$
		\STATE $z[0:n-1]\leftarrow \text{VecSumErrBranch}(y[0:m-1],n)$
		\RETURN  $z_0,z_1,\cdots,z_{n-1}$   
	\end{algorithmic}
\end{algorithm}

If the relationship between the elements of $\{x_i\}_{i=0}^{m-1}$ cannot be established in advance, we can invoke the third condition for VecSum and apply VecSum to $(x_0,x_1)$, then to $(x_0,x_1,x_2)$, appending one more term $x_i$ at the tail each time, until all $x_i$ have been added. This yields Algorithm \ref{algo:renrmrand}, which is suitable for random inputs.

\begin{algorithm}[H]
	\algsetup{linenosize=\small}
	\small
	\caption{RenormalizeRandom$(x_0,x_1,\cdots,x_{m-1},n)$\cite{Popescu2017Phd}}\label{algo:renrmrand}
	\begin{algorithmic}[1]
		\STATE $e_0^{(0)}\leftarrow x_0$
		\FOR {$i\leftarrow 1$ \TO ${m-1}$}
		\STATE $e^{(i)}[0:i]\leftarrow \mathrm{VecSum}(e^{(i-1)}[0:i-1],x_i)$
		\ENDFOR
		\STATE $r[0:n-1]\leftarrow \text{VecSumErrBranch}(e^{(m-1)}[0:m-1],n)$
		\RETURN  $r_0,r_1,\cdots,r_{n-1}$   
	\end{algorithmic}
\end{algorithm}

Fabiano et al. \cite{JMuller2019TripleDouble} showed that, under certain conditions, the output of Algorithm \ref{algo:renrm} is $\mathcal{P}$-nonoverlapping, i.e., $|r_i|<\mathrm{ulp}(r_{i-1})$, which is tighter than ulp-nonoverlapping. 

Having introduced the renormalization algorithms, we next describe addition algorithms. The core idea of these arithmetic algorithms is as follows: transform the result into an array, either error-free or with bounded error, that is as short as possible and as close as possible to the conditions required by VecSum; then apply a renormalization algorithm to the array. 

We assume that the inputs of these arithmetic algorithms are ulp-nonoverlapping and that all input and output expansions have the same length. Input expansions cannot contain interleaving zeros; zeros may appear only at the end.

To build an addition algorithm, note the last condition for VecSum: we can merge $\{x_i\}$ and $\{y_i\}$ and then sort the result by absolute value. Since $\{x_i\}$ and $\{y_i\}$ are already ulp-nonoverlapping, the merge sort algorithm (Algorithm \ref{algo:merge}) suffices. After that, Algorithm \ref{algo:renrm} is called to make the returned array ulp-nonoverlapping and truncate it to the desired length. This gives the merge-based Algorithm \ref{algo:mergeadd}. Since $\{r_i\}_{i=0}^{n-1}$ are the first $n$ terms of the ulp-nonoverlapping expansion $\{r_i\}_{i=0}^{2n-1}$ satisfying $\sum_{i=0}^{2n-1}r_i=\sum_{i=0}^{2n-1}f_i=x+y$, the relative error of Algorithm \ref{algo:mergeadd} arises only from the truncation step and is bounded by $\frac {\sum_{i=n}^\infty 2^{-i(p-1)}} { 1 - \sum_{i=1}^\infty 2^{-i(p-1)}}\leq \frac 1 {1-2^{2-p}} 2^{-n(p-1)}$ (using $|r_{i}|\leq \mathrm{ulp}(r_{i-1}) \leq 2^{-(p-1)}|r_{i-1}|\leq 2^{-i(p-1)}|r_0|$ for $i>0$).    

\begin{algorithm}[H]
	\algsetup{linenosize=\small}
	\small
	\caption{Merge$(x_0,\cdots,x_{n-1},y_0,\cdots,y_{n-1})$}\label{algo:merge}
	\begin{algorithmic}[1]
		\STATE $i,j\leftarrow 0$
		\FOR { $t\leftarrow  0$ \TO $2n-1$  }
		\IF {$i=n$ \OR ($j<n$ \AND $|y_j|>|x_i|$) }
		\STATE $f_t\leftarrow y_j;\ j\leftarrow j+1$
		\ELSE
		\STATE $f_t\leftarrow x_i;\ i\leftarrow i+1$
		\ENDIF
		\ENDFOR
		\RETURN  $f_0,f_1,\cdots,f_{2n-1}$   
	\end{algorithmic}
\end{algorithm}

\begin{algorithm}[H]
	\algsetup{linenosize=\small}
	\small
	\caption{MergeAdd$(x_0,\cdots x_{n-1},y_0,\cdots, y_{n-1})$}\label{algo:mergeadd}
	\begin{algorithmic}[1]
		\STATE $f[0:2n-1]\leftarrow \text{Merge}(x[0:n-1],y[0:n-1])$
		\STATE $r[0:n-1] \leftarrow \text{Renormalize}(f[0:2n-1],n)$
		\RETURN $r_0,r_1,\cdots,r_{n-1}$
	\end{algorithmic}
\end{algorithm}

Another addition algorithm, Algorithm \ref{algo:QDaddRandom}, is a generalization of the quad-double addition algorithm implemented in the QD library \cite{Bailey2000QD}. It has the absolute error bound $|(x+y+\epsilon_1(|x_0|+|y_0|))(1+\epsilon_2)-(x+y)|$, where $|\epsilon_1|\leq \frac {n-1} {2} 2^{-(n+1)(p-1)}$ is the sum of the roundoff errors in step 8 \cite{Popescu2017Phd}, and $|\epsilon_2|\leq \frac 1 {1-2^{2-p}} 2^{-n(p-1)}$ is the relative error of the truncation step in renormalization. Therefore, the relative error of Algorithm \ref{algo:QDaddRandom} is $\left |\epsilon_1(1+\epsilon_2)\frac {|x_0|+|y_0|} {x+y} +\epsilon_2\right|$. Although the relative error depends on the amount of cancellation, measured by $\frac {|x_0|+|y_0|} {|x+y|} $, it remains of the order of $2^{-n(p-1)}$ as long as $\frac {|x_0|+|y_0|} {|x+y|} <2^{p-1} $. 

At first sight, Algorithm \ref{algo:QDaddRandom} performs more floating-point operations than Algorithm \ref{algo:mergeadd}. However, Algorithm \ref{algo:QDadd} is faster in most practical cases, owing to better instruction-level parallelism. 

The original addition function in the QD library, named \texttt{qd\_real::sloppy\_add}, implements Algorithm \ref{algo:QDadd}, not Algorithm \ref{algo:QDaddRandom}. The difference between the two algorithms lies in the renormalization step. Popescu \cite{Popescu2017Phd} reports being unable to prove a strict relationship between the $f_i$, and therefore chooses Algorithm \ref{algo:QDaddRandom}. Here we show that $\mathrm{uls}(f_i)\geq \mathrm{ulp}(f_{i+1})$, so it is safe to use Algorithm \ref{algo:QDadd}. Moreover, Fast2Sum can even be used instead of 2Sum inside VecSum in Algorithm \ref{algo:renrm}.

\begin{algorithm}[H]
	\algsetup{linenosize=\small}
	\small
	\caption{QDAdd$(x_0,\cdots,x_{n-1},y_0,\cdots,y_{n-1})$}\label{algo:QDadd}
	\begin{algorithmic}[1]
		\STATE $(f_0,e_0)\leftarrow \mathrm{2Sum}(x_0,y_0)$
		\FOR {$i\leftarrow 1$ \TO $n-1$}
		\STATE $(s_i,e_i)\leftarrow \mathrm{2Sum}(x_i,y_i)$
		\STATE $(f_i,e[0:i-1])\leftarrow \mathrm{VecSum}(s_i,e[0:i-1])$
		\ENDFOR
		\STATE $f_n \leftarrow 0$
		\FOR {$i\leftarrow 0$ \TO $n-1$}
		\STATE $f_n \leftarrow f_n + e_i$
		\ENDFOR
		\STATE $r[0:n-1]\leftarrow\text{Renormalize}(f[0:n],n)$
		\RETURN $r_0,r_1,\cdots,r_{n-1}$
	\end{algorithmic}
\end{algorithm}

\begin{algorithm}[H]
	\algsetup{linenosize=\small}
	\small
	\caption{QDAddRandom$(x_0,\cdots,x_{n-1},y_0,\cdots,y_{n-1})$}\label{algo:QDaddRandom}
	\begin{algorithmic}[1]
		\STATE $(f_0,e_0)\leftarrow \mathrm{2Sum}(x_0,y_0)$
		\FOR {$i\leftarrow 1$ \TO $n-1$}
		\STATE $(s_i,e_i)\leftarrow \mathrm{2Sum}(x_i,y_i)$
		\STATE $(f_i,e[0:i-1])\leftarrow \mathrm{VecSum}(s_i,e[0:i-1])$
		\ENDFOR
		\STATE $f_n \leftarrow 0$
		\FOR {$i\leftarrow 0$ \TO $n-1$}
		\STATE $f_n \leftarrow f_n + e_i$
		\ENDFOR
		\STATE $r[0:n-1]\leftarrow\text{RenormalizeRandom}(f[0:n],n)$
		\RETURN $r_0,r_1,\cdots,r_{n-1}$
	\end{algorithmic}
\end{algorithm}

To analyze the VecSum algorithm, we first examine the first term of the output, i.e., $e_0$. It is easy to see that $e_0=\mathrm{RN}(x_0+\mathrm{RN}(x_1+\cdots+\mathrm{RN}(x_{m-2}+x_{m-1})\cdots))$, which equals the result of summing the array $\{x_i\}_{i=0}^{m-1}$ from the tail $x_{m-1}$ to the head $x_0$. We denote this by RORS (Reverse Ordinary Recursive Summation), and by ORS (Ordinary Recursive Summation) when the summation starts from the head.
\begin{algorithm}[H]
	\algsetup{linenosize=\small}
	\small
	\caption{RORS$(x_0,\cdots,x_{n-1})$}\label{algo:RORS}
	\begin{algorithmic}[1]
		\STATE $r\leftarrow x_{n-1}$
		\FOR {$i\leftarrow n-2$ \TO $0$}
		\STATE $r\leftarrow \mathrm{RN}(r+x_{i})$
		\ENDFOR
		\RETURN $r$
	\end{algorithmic}
\end{algorithm}

\begin{algorithm}[H]
	\algsetup{linenosize=\small}
	\small
	\caption{ORS$(x_0,\cdots,x_{n-1})$}\label{algo:ORS}
	\begin{algorithmic}[1]
		\STATE $r\leftarrow x_0$
		\FOR {$i\leftarrow 1$ \TO $n-1$}
		\STATE $r\leftarrow \mathrm{RN}(r+x_{i})$
		\ENDFOR
		\RETURN $r$
	\end{algorithmic}
\end{algorithm}

To derive error bounds for ORS and RORS, we introduce the $\theta_n$ and $\gamma_n$ notation defined by Higham \cite{Higham2002Accuracy2nd}.
\begin{definition}
	For $\varepsilon_i$ such that $\varepsilon_{i}\leq \mathbf{u}$, $0\leq i<n$, and assuming $n\mathbf{u}<1$,
	\begin{equation*}
	 \prod_{i=0}^{n-1} (1+\varepsilon_i)^{\pm1} =1+\theta_n,
	\end{equation*} 
	where
	\begin{equation*}
	|\theta_n|\leq \frac {n\mathbf{u}} {1-n\mathbf{u}}.
	\end{equation*}
	And $\gamma_n$ is defined as
	\begin{equation*}
	\gamma_n \triangleq \frac {n\mathbf{u}} {1-n\mathbf{u}}.
	\end{equation*}
\end{definition}

Both ORS and RORS satisfy the same forward error bound, $\gamma_{n-1}\sum_{i=0}^{n-1}|x_i|$ \cite{HandbookFP2018}.
\begin{equation}
\left| \mathrm{ORS}(\{x_i\}_{i=0}^{n-1})-\sum_{i=0}^{n-1}x_i\right|\leq\gamma_{n-1}\sum_{i=0}^{n-1}|x_i|,
\end{equation}
\begin{equation}
\left| \mathrm{RORS}(\{x_i\}_{i=0}^{n-1})-\sum_{i=0}^{n-1}x_i\right|\leq\gamma_{n-1}\sum_{i=0}^{n-1}|x_i|.
\end{equation}

\section{Review of Renormalization Algorithms}

This section mainly discusses the behavior of VecSum and the F-nonoverlapping and ulp-nonoverlapping conditions.

Since zeros in the input FPE do not affect the VecSum result, unless stated otherwise, all subsequent proofs assume that the input FPE contains no zeros; intermediate results and the output FPE may contain zeros.

In addition, we restrict the length of the FPE to $n\leq p$, which almost always holds in practice. We also require $n\geq 3$.

\begin{lemma}\label{lem:adduls}
	Let $a,b\in\mathbb{F}_p$ and $(s,t)\leftarrow \mathrm{2Sum}(a,b)$, and let $c=\min(\mathrm{uls}(a),\mathrm{uls}(b))$. Then either $\mathrm{uls}(s)\geq c$ or $s=0$, and either $\mathrm{uls}(t)\geq c$ or $t=0$.
\end{lemma}

\begin{lemma}
If $\{x_i\}$ is F-nonoverlapping, then $|x_{i+1}|\leq \frac  1 2 \mathrm{uls}(x_{i})\leq \frac {|x_i|} 2$.
\end{lemma}

The proofs of the two lemmas above are elementary and are therefore omitted.

We first analyze the properties of the VecSum output.
\begin{lemma}[Sum of an F-nonoverlapping expansion]
	For an F-nonoverlapping expansion $\{x_i\}_{i=0}^{n-1}$ with $|x_1|\leq \frac 1 2 \mathrm{ulp}(x_0)$, suppose $n\leq p+1$. If $\{y_i\}_{i=1}^{n-1} \leftarrow \mathrm{VecSum}(\{x_i\}_{i=1}^{n-1})$, then $|y_1|<\mathrm{ulp}(x_0)$.
\end{lemma}
\begin{proof}
	By induction, $|x_i|\leq 2^{-i} \mathrm{ulp}(x_0)$, which follows from the monotonicity of RN. 
	\begin{align*}
	|y_1| &= |\mathrm{RN}(x_1 + \mathrm{RN}(x_2 + \cdots +\mathrm{RN}(x_{n-2}+x_{n-1})\cdots)| \\
	&\leq |\mathrm{RN}(|x_1| + \mathrm{RN}(|x_2| + \cdots +\mathrm{RN}(|x_{n-2}|+|x_{n-1}|)\cdots))| \\
	&\leq \mathrm{ulp}(x_0)|\mathrm{RN}(2^{-1} + \mathrm{RN}(2^{-2} + \cdots +\mathrm{RN}(2^{-(n-2)}+2^{-(n-1)})\cdots))|\\
	& = \mathrm{ulp}(x_0)(1-2^{-(n-1)})\le \mathrm{ulp}(x_0)(1-2^{-p}) \\
	& < \mathrm{ulp}(x_0).
	\end{align*}
\end{proof}

Assume that the FPE $\{x_i\}_{i=0}^{n-1}$ satisfies $|x_{i+1}|\leq q|x_i|$ with $0<q<\frac 1 {1+\mathbf{u}}$. We now estimate its RORS result.
	\begin{align}
	|\mathrm{RORS}(\{x_i\}_{i=1}^{n})| &= |\mathrm{RN}(x_1 + \mathrm{RN}(x_2 + \cdots +\mathrm{RN}(x_{n-2}+x_{n-1})\cdots)|   \notag \\
	&\leq |\mathrm{RN}(|x_1| + \mathrm{RN}(|x_2| + \cdots +\mathrm{RN}(|x_{n-2}|+|x_{n-1}|)\cdots)| \notag \\
	&\leq\mathrm{RN}( |x_0|q + \mathrm{RN}( |x_0|q^2 + \cdots +\mathrm{RN}( |x_0|q^{n-2}+ |x_0|q^{n-1})\cdots))\label{eq:sumqfpe}\\
	&\leq |x_0|(1+\mathbf{u})(q+(1+\mathbf{u})  (q^2 + \cdots +(1+\mathbf{u})(q^{n-2}+q^{n-1})\cdots))\notag\\
	& = |x_0| ( (1+\mathbf{u})q + (1+\mathbf{u})^2q^2 + \cdots + (1+\mathbf{u})^{n-2}q^{n-2}+(1+\mathbf{u})^{n-2}q^{n-1} )\notag \\
	& < |x_0| \Sigma_{i=1}^{n-1} (1+\mathbf{u})^iq^i \notag\\
	& = |x_0| \frac {(1+\mathbf{u})q- (1+\mathbf{u})^nq^n }  {1-(1+\mathbf{u})q } \notag\\
	& = |x_0| \left (1-\frac 1 {1-(1+\mathbf{u})q} ( 1-2(1+\mathbf{u})q+(1+\mathbf{u})^nq^n ) \right ).\notag
	\end{align}
	Since $ (1+\mathbf{u})q<1$, the inequality $|\mathrm{RORS}(\{x_i\}_{i=1}^{n})|\leq |x_0|$ is  equivalent to $ 1-2(1+\mathbf{u})q+(1+\mathbf{u})^nq^n  \geq 0$. 
	
	\begin{lemma}\label{lem:f-nonoverlap-rors}
		If $\{x_i\}_{i=0}^{n-1}$ is F-nonoverlapping and $n\leq p$, then $|\mathrm{RORS}(\{x_i\}_{i=1}^{n})|< \mathrm{uls}(x_0)$.
	\end{lemma}
\begin{proof}
		If $\{x_i\}_{i=0}^{n-1}$ is F-nonoverlapping, substituting $q =\frac  1 2$ gives $ (1+\mathbf{u})^n 2^{-n}\geq2^{-n}\geq2^{-p}=\mathbf{u} $ and $1-2(1+\mathbf{u})q+(1+\mathbf{u})^nq^n =(1+\mathbf{u})^nq^n-\mathbf{u}\geq 0$, hence $|\mathrm{RORS}(\{x_i\}_{i=1}^{n})|< |x_0|$. This bound can be refined to $|\mathrm{RORS}(\{x_i\}_{i=1}^{n})|< \mathrm{uls}(x_0)$ by replacing $x_0$ with $\mathrm{uls}(x_0)$ in \eqref{eq:sumqfpe}; that is, when $\{e_i\}\leftarrow\mathrm{VecSum}(\{x_i\})$, $|x_{i-1}|>|t_i|$. A simpler proof is by induction: $|\mathrm{RORS}(x_{i+1},x_{i+2},\cdots)|=|\mathrm{RN}(x_{i+1}+\mathrm{RORS}(x_{i+2}+\cdots))|< \mathrm{RN}(2^{-1}\mathrm{uls}(x_i),2^{-1}\mathrm{uls}(x_{i+1}))<\mathrm{uls}(x_i) $.
\end{proof}	
	If $\{x_i\}_{i=0}^{n-1}$ overlaps by at most $d\leq p-2$ bits, then $q =\frac {2^p-1} {2^{p+1}}  = \frac 1 2 (1-\mathbf{u})<1/2$, and we likewise have $|\mathrm{RORS}(\{x_i\}_{i=1}^{n})|<|x_0|$.

\begin{lemma}\label{lem:Fbasecond}
	If $\{x_i\}_{i=0}^{n-1}$ satisfies $\mathrm{uls}(x_j)\geq \mathrm{ulp}(\mathrm{RORS}(\{x_i\}_{i=j+1}^{n-1}))$ for $j=0,1,\cdots,n-2$, then $\{y_i\}_{i=0}^{n-1}\leftarrow\mathrm{VecSum}(\{x_i\}_{i=0}^{n-1})$ is F-nonoverlapping, and the 2Sum in VecSum can be replaced by Fast2Sum.
\end{lemma}
\begin{proof}
Since $t_j=\mathrm{RORS}(\{x_i\}_{i=j}^{n-1})$ in the VecSum algorithm (Algorithm \ref{algo:vecsum}), i.e., $\mathrm{uls}(x_j)>\mathrm{ulp}(t_{j+1})$, it follows from Lemma \ref{lem:fast2sumextra} that $\mathrm{2Sum}(x_i,t_{i+1})$ can be replaced by $\mathrm{Fast2Sum}(x_i,t_{i+1})$. 
	For $y_{i+1}\neq0$, $|y_{i+1}|\leq \frac 1 2 \mathrm{ulp}(t_i)$.
	$(t_{i-1},y_i)\leftarrow \mathrm{2Sum}(x_{i-1},t_i)$. By the condition, $\mathrm{uls}(x_{i-1})\geq \mathrm{ulp}(t_i)$.
	By Lemma \ref{lem:adduls}, $\mathrm{uls}(y_i)\geq \min(\mathrm{uls}(x_{i-1}),\mathrm{uls}(t_i)) \geq \mathrm{ulp}(t_i) \geq 2 |y_{i+1}|$, or $y_i=0$.
	When $y_i\neq 0$, $|y_{i+1}|\leq \frac 1 2 \mathrm{uls}(y_i)$, so the VecSum output $\{y_i\}_{i=0}^{n-1}$ is F-nonoverlapping; the same holds when $y_{i+1}=0$ or $y_i=0$.  
\end{proof}
This lemma gives a sufficient condition for the VecSum output to be F-nonoverlapping.

\begin{theorem}\label{thm:vecsum_cond_random}%[condition.c]
	Suppose $(x_0,x_1,\cdots,x_{n-1})$ is F-nonoverlapping and $s$ is a random number. Then $(y_0,y_1,\cdots,y_{n-1},y_n)\leftarrow \mathrm{VecSum}(x_0,x_1,\cdots,x_{n-1},s)$ is F-nonoverlapping.
\end{theorem}
\begin{proof}
	The first step in VecSum is $(t_{n-1},y_n)\leftarrow\mathrm{2Sum}(x_{n-1},s)$, and $t_{n-1}$ is the floating-point number closest to $x_{n-1}+s$, \[
|y_n|=| (x_{n-1}+s)-t_{n-1} |\leq |(x_{n-1}+s)-s | =|x_{n-1}|.
\]
	$|y_n|\leq \frac 1 2 \mathrm{ulp}(t_{n-1})\leq \frac 1 2 \mathrm{uls}(t_{n-1})$; moreover,
	$|y_n|\leq |x_{n-1}| \leq \frac 1 2 \mathrm{uls}(x_{n-2})$. Since $(t_{n-2},y_{n-1})\leftarrow\mathrm{2Sum}(x_{n-2},t_{n-1})$, Lemma \ref{lem:adduls} gives
	\begin{align*}
	\mathrm{uls}(y_{n-1})
	&\geq \min(\mathrm{uls}(t_{n-1}),\mathrm{uls}(x_{n-2})) \\
	&\geq 2|y_n|
	\end{align*}
or $y_{n-1} = 0$. When $y_{n-1}=0$: if all $y_i=0$ for $0\leq i<n-1$, then $(y_0,y_1,\cdots,y_{n-1},y_n)$ still satisfies the F-nonoverlapping condition; otherwise, let $j=\arg\min\{i\mid y_i\neq0,\,0\leq i<n-1\}$; Lemma \ref{lem:adduls} again yields $\mathrm{uls}(y_i)\geq\min(\mathrm{uls}(t_{n-1}),\mathrm{uls}(x_{n-2}),\mathrm{uls}(x_{n-3}),\cdots,\mathrm{uls}(x_{i-1}) )\geq 2|y_n|$. By induction, $(y_0,y_1,\cdots,y_{n-1},y_n)$ is F-nonoverlapping. The case $y_n=0$ does not affect this conclusion.
\end{proof}	 
	
This theorem also implies that if $\{x_i\}$ is F-nonoverlapping, then $\mathrm{VecSum}(\{x_i\})$ and every one of its subsets are F-nonoverlapping.
	
Under this condition, one cannot guarantee that there is no $0 <  i <n$ such that $y_{i-1}y_{i}<0$ and $|y_{i-1}|=2^{k+1},|y_{i}|=2^{k}$ with $k\in\mathbb{Z}$.
	For example, $\{x_i\}=(-2,3\cdot 2^{1-p},-2^{-p}, 2-2^{2-p})$ yields the output $\{y_i\}=(0,0,2^{1-p},-2^{-p})$.

\begin{theorem}\label{thm:vec-on-F-toF}%[condition.e]
	Suppose $(x_0,x_1,\cdots,x_{n-1})$ is F-nonoverlapping, $n\leq p$, and ties are broken by the ties-to-even rule. Then $(y_0,y_1,\cdots,y_{n-1})\leftarrow \mathrm{VecSum}(x_0,x_1,\cdots,x_{n-1})$ is F-nonoverlapping, there is no $0 \leq j < i <n$ such that $ |y_j|=2^{k+1}$, $y_m=0$ for $j<m<i$, and $|y_{i}|=2^{k}$ with $k\in\mathbb{Z}$, and Fast2Sum can be used in VecSum. 
\end{theorem}	
\begin{proof}
	By Theorem \ref{thm:vecsum_cond_random}, $\{y_i\}_{i=0}^{n-1}$ is F-nonoverlapping. By Lemma \ref{lem:f-nonoverlap-rors}, $|t_i|=\mathrm{RORS}(x_{i+1},\cdots,x_{n-1})<\mathrm{uls}(x_{i-1})$, so Fast2Sum can be safely used here by Lemma \ref{lem:fast2sumextra}. 
	
	It is easy to obtain $|y_i|\leq \frac 1 2 \mathrm{ulp}(t_{i-1})$ and $|t_{i-1}| <\mathrm{uls}(x_{i-2})$. If $|y_i|=\frac 1 2 \mathrm{ulp}(t_{i-1})$, then, by the ties-to-even rule, the last significant digit of $t_{i-1}$ must be zero; in other words, $\mathrm{uls}(t_{i-1})\geq 2\mathrm{ulp}(t_{i-1})=4|y_i|$. Moreover, $ \mathrm{uls}(x_j)>|t_{i-1}|\geq \mathrm{uls}(t_{i-1})\geq 4|y_i|$ for $j <i-1$. Thus, by Lemma \ref{lem:adduls}, we obtain $|y_j|\geq 4|y_i|$ or $y_j=0$ for $0\leq j < i$. Otherwise, if $|y_i|=2^k$ with $k\in\mathbb{Z}$, then necessarily $|y_i|\leq\frac 1 4\mathrm{ulp}(t_{i-1})$, and again we get $|y_j|\geq 4|y_i|$ or $y_j=0$ for $0\leq j < i$. Hence there is no $0 <  i \leq n-1$ such that $ |y_{i-1}|=2^{k+1}$ and $|y_{i}|=2^{k}$ with $k\in\mathbb{Z}$. If $y_{i-1}=0,\cdots$ precede $y_i$, the same argument shows that there is no $0 \leq j < i <n$ such that $ |y_j|=2^{k+1}$, $y_m=0$ for $j<m<i$, and $|y_{i}|=2^{k}$ with $k\in\mathbb{Z}$.
\end{proof}
This theorem shows that the output remains F-nonoverlapping after repeated applications of VecSum.

We now prove that, for an F-nonoverlapping FPE, repeated applications of VecSum guarantee an output that is $\mathcal{P}$-nonoverlapping or ulp-nonoverlapping (Algorithm \ref{algo:multivecsum}).  

\begin{theorem}
	Assume that $\{h_i\}_{i=0}^{n-1}$ is F-nonoverlapping, $n\geq 3$, and either $h_0h_1>0$ or there is no $k\in\mathbb{Z}$ such that $|h_0|=2^k$ and $|h_1|=2^{k-1}$. Let $ \{f_i\}_{i=0}^{n-1}\leftarrow \mathrm{VecSum}(\{h_i\}_{i=0}^{n-1})$ and $\{g_i\}_{i=1}^{n-1}\leftarrow\mathrm{VecSum}(\{f_i\}_{i=1}^{n-1})$. Then $|g_1| < \mathrm{ulp}(f_0)$.
\end{theorem}
\begin{proof}
	In $\{g_i\}_{i=1}^{n-1}\leftarrow\mathrm{VecSum}(\{f_i\}_{i=1}^{n-1})$, denote by $\hat{t}_{i+1}$ the quantity $t_i$ of Algorithm \ref{algo:vecsum}, to distinguish it from the $t_i$ in $ \{f_i\}_{i=0}^{n-1}\leftarrow \mathrm{VecSum}(\{h_i\}_{i=0}^{n-1})$.
	By Theorem \ref{thm:vec-on-F-toF}, both $ \{f_i\}_{i=0}^{n-1}$ and $\{g_i\}_{i=1}^{n-1}$ are F-nonoverlapping.
	% Since $|t_1| < |e_0|$. 
	
	Clearly $f_0=\mathrm{RORS}(\{h_i\}_{i=0}^{n-1})=\mathrm{RN}(h_0 +\mathrm{RORS}(\{h_i\}_{i=1}^{n-1}) )\neq0$.
	
	If $f_1\neq 0$, then $|f_1|\leq \frac 1 2 \mathrm{ulp}(f_0)$ and $\hat{t}_2=\mathrm{RORS}(f_2,f_3,\cdots,f_{n-1})$. Since $\{f_i\}$ is F-nonoverlapping, Lemma \ref{lem:f-nonoverlap-rors} gives $|\hat{t}_2|<|f_1|\leq \frac 1 2 \mathrm{ulp}(f_0)$. Thus $|g_1|=\mathrm{RN}(f_1+\hat{t}_2)\leq\mathrm{ulp}(f_0)$. Moreover,
	\begin{align*}
	\hat{t}_2&=\mathrm{RORS}(f_2,f_3,\cdots,f_{n-1})\\
	& \leq \mathrm{RN}(\mathrm{uls}(f_1)(2^{-1}+2^{-2}+\cdots+2^{2-n})) \\
	& \leq \mathrm{RN}(\mathrm{uls}(f_1)(2^{-1}+2^{-2}+\cdots+2^{2-p}))\\ 
	& \leq \mathrm{uls}(f_1)(2^{-1}+2^{-2}+\cdots+2^{2-p})\\
	& \leq \frac 1 2 \mathrm{ulp}(f_0)(2^{-1}+2^{-2}+\cdots+2^{2-p}).
	\end{align*}
	
	\begin{align*}
	|g_1|&=\mathrm{RN}(f_1+\hat{t}_2)\\
	 & \leq \mathrm{RN}\left(\frac 1 2 \mathrm{ulp}(f_0) +\frac 1 2 \mathrm{ulp}(f_0)(2^{-1}+2^{-2}+\cdots+2^{2-p}) \right)\\
	 & = \frac 1 2 \mathrm{ulp}(f_0) \mathrm{RN}(1+2^{-1}+2^{-2}+\cdots+2^{2-p}) \\
	 & = \frac 1 2 \mathrm{ulp}(f_0)(1+2^{-1}+2^{-2}+\cdots+2^{2-p}) \\
	 & <  \mathrm{ulp}(f_0).
	\end{align*}
	
	If $f_1=0$: if all $f_i$ with $i>1$ are zero, then clearly $|g_1|=0<\mathrm{ulp}(f_0)$. Since $\{f_i\}_{i=2}^{n-1}$ with $k>0$ zeros is equivalent to $\{\tilde{f}_i\}_{i=2}^{n-1-k}$ after deleting all the zeros, we may assume $f_i\neq 0$ for $i>1$. Here $f_i$ is the roundoff error of $\mathrm{RN}(h_{i-1}+t_i)$. From Lemma \ref{lem:f-nonoverlap-rors} we know that $|t_{i}|<\mathrm{uls}(h_{i-1})$. If $h_{i-1}=2^k$ with $k\in\mathbb{Z}$, then $\mathrm{ufp}(t_{i-1})=\mathrm{ufp}(\mathrm{RN}(h_{i-1}+t_i))\leq\mathrm{ufp}(h_{i-1})$; the proof is analogous to the estimate of $|g_1|$ above and is omitted, and we obtain $|f_i| \leq \frac 1 2 \mathrm{ulp}(t_{i-1}) \leq \frac 1 4 \mathrm{ulp}(f_{i-2})$. Otherwise, $\mathrm{ufp}(\mathrm{RN}(f_{i-1}+t_i))\leq 2\mathrm{ufp}(f_{i-1})$, so $|f_i| \leq  \mathrm{ulp}(f_{i-1}) \leq \frac 1 4 \mathrm{ulp}(f_{i-2})$. In summary, $|f_i|\leq \frac 1 4 \mathrm{ulp}(f_{i-2})$.

	\begin{align*}
	|g_1|&=|\mathrm{RORS}(f_1,f_2,f_3,\cdots,f_{n-1})| \\
	&=|\mathrm{RORS}(f_2,f_3,\cdots,f_{n-1})| \\
	&\leq  \frac 1 4 \mathrm{RORS}(\mathrm{ulp}(e_0),\mathrm{ulp}(2^{-1}\cdot {e_0} ),\cdots,\mathrm{ulp}(2^{-(n-2)}\cdot{e_0} )  ) \\
	& =  \frac 1 4\mathrm{ulp}(e_0) \mathrm{RORS}(1, 2^{-1},\cdots,2^{-(n-2)}) \\
	& \leq \frac 1 2 \mathrm{ulp}(e_0).
	\end{align*}
	\begin{itemize}	
		\item If $|e_0|\neq 2^k$ for all $k\in\mathbb{Z}$: since $|t_1|< \mathrm{uls}(e_0)$, $\mathrm{ulp}(f_0)\geq \mathrm{ulp}(|e_0|-|t_1|)=\mathrm{ulp}(e_0)$, and hence $|g_1| < \frac 1 2 \mathrm{ulp}(e_0) \leq \mathrm{ulp}(g_0)$.
		\item Otherwise $|e_0|=2^k$. If $e_0t_1\geq 0$, then obviously $|g_1| < \frac 1 2 \mathrm{ulp}(e_0) \leq \mathrm{ulp}(g_0)$. 
		
		We now consider the case $|e_0|=2^k$ and $e_0t_1<0$.

		\begin{itemize}
			\item If $|e_1|< 2^{k-1}$: then $|e_1|<\frac 1 2 |e_0|$, $|t_2|< \mathrm{uls}(e_1)$, and $|t_1|\leq 2\mathrm{ufp}(e_1) \leq \frac 1 2 |e_0|$, so $|g_0|=|\mathrm{RN}(e_0+t_1)|\geq \frac 1 2 |e_0|$ and $|g_1| < \frac 1 2 \mathrm{ulp}(e_0) \leq \mathrm{ulp}(g_0)$.
			
			\item If $|e_1|=2^{k-1}$: similarly, if $|e_2| < 2^{k-2}$, then $|t_2|< \frac {|e_1|} 2\leq \frac 1 4 |e_0|$ and $|g_0| = |\mathrm{RN}(e_0+\mathrm{RN}(e_1+t_2))|\geq \frac 1 4 |e_0|$. Now if $f_2=0$, then $|g_1|=|\mathrm{RORS}(f_3,f_4,\cdots,f_{n-1})| < \frac 1 4 \mathrm{ulp}(e_0)  \leq \mathrm{ulp}(g_0)$.
			% F-nonoverlap can't ensure $\{g_i\}$ is ulp-nonoverlapping. eg., $\{-1, 2^{-1}, 2^{-2}, 2^{-53}+2^{-54} \} \rightarrow \{-2^{-2}+2^{-52}  , -2^{-54},0,0\}$. 

			If $f_2\neq 0$, we can give a counterexample: $\{e_i\}=(-2^3,2^2,1+2^{-(p-2)}, 2^{-(p+1)})$ gives $\{g_i\}=(-3, 2^{-(p-2)}+2^{-(p+1)},0,0)$ and $\mathrm{ulp}(3)=2^{-(p-2)}$.

		\end{itemize}
	\end{itemize}	
	
\end{proof}

\begin{theorem}
	The minimum length of an F-nonoverlapping FPE $\{h_i\}_{i=0}^{n-1}$ such that, with $ \{f_i\}_{i=0}^{n-1}\leftarrow \mathrm{VecSum}(\{h_i\}_{i=0}^{n-1})$ and $\{g_i\}_{i=1}^{n-1}\leftarrow\mathrm{VecSum}(\{f_i\}_{i=1}^{n-1})$, one has $|g_1| > \mathrm{ulp}(f_0)$, is $n=4$.
\end{theorem}
\begin{proof}
	A counterexample for $n=4$ is $\{h_i\}=(-2^3,2^2,1+2^{-(p-2)}, 2^{-(p+1)})$, for which $\{g_i\}=(-3, 2^{-(p-2)}+2^{-(p+1)},0,0)$ and $\mathrm{ulp}(3)=2^{-(p-2)}$.
	Hence the minimum length is at most 4. For $n=1,2$ it is obvious that $|g_1|\leq \mathrm{ulp}(f_0)$, so it remains to prove that $|g_1|\leq \mathrm{ulp}(f_0)$ when $n=3$. 
	
	By the above theorem, the only possible case is $\{h_0=4c,h_1=-2c,h_2\}$ with $c=\pm 2^k$, $k\in\mathbb{Z}$, and $|h_2|\leq c$. 
	\begin{itemize}
		\item $h_1h_2\geq 0$. $|t_1|\leq 3|c|$, $|f_0|\geq |c|$. By Sterbenz's lemma, $|t_1|>|h_0|/2$ and $|f_1|=0$. Also $|f_2|\leq \frac 1 2\mathrm{ulp}(t_1)\leq \mathrm{ulp}(c)$, so $|g_1|=|\mathrm{RN}(f_1+f_2)| =|f_2|\leq \mathrm{ulp}(c) \leq \mathrm{ulp}(f_0)$. 
		\item $h_1h_2<0$. $|t_1|\leq 2|c|$, $|f_0|\geq 2|c|$. $|f_2|\leq \frac 1 2 \mathrm{ulp}(t_1) \leq \frac 1 2 \mathrm{ulp}(2c) \leq \frac 1 2 \mathrm{ulp}(f_0)$, and $|f_1|\leq \frac 1 2 \mathrm{ulp}(f_0)$.  $|g_1|=|\mathrm{RN}(f_1+f_2)| \leq \mathrm{ulp}(f_0)$. 
	\end{itemize}
\end{proof}

\begin{remark}
	If the F-nonoverlapping FPE  $\{h_i\}_{i=0}^{n-1}$  is the result of VecSum, then its first two terms $h_0,h_1$ cannot satisfy $|h_1|=|h_0|/2$. 
\end{remark}

\begin{theorem}
	Assume that $\{x_i\}_{i=0}^{n-1}$ is F-nonoverlapping and $n\leq p$. If $\{y_i\}\leftarrow \text{VecSumErrBranch}(\{x_i\}_{i=0}^{n-1})$, then $\{y_i\}$ is $\mathcal{P}$-nonoverlapping.
\end{theorem}
\begin{proof}
	Suppose $(y_{i-1},t_i)\leftarrow \mathrm{Fast2Sum}(t_{i-1},x_{i})$ with $t_i\neq 0$. 
	
	Since $|y_i| \leq \mathrm{ORS}(|t_i|,|x_{i+1}|,\cdots,|x_{n-1}|)$, to prove $\mathrm{ulp}(y_{i-1})> |y_{i}|$ it suffices to show that $\mathrm{ORS}(|t_i|,|x_{i+1}|,\cdots,|x_{n-1}|) < \mathrm{ulp}(y_{i-1})$.
	
	$|t_i| \geq \mathrm{uls}(t_i) \geq \min_{j=0}^{i}(\mathrm{uls}(x_j)) \geq 2 |x_{i+1}|$. In other words, $(t_i,x_{i+1}, x_{i+2},\cdots,x_{n-1})$ is F-nonoverlapping, and because $|t_i|\leq \frac 1 2 \mathrm{ulp}(y_{i-1})$, we get $|f_i| \leq \mathrm{ORS}(|t_i|,|x_{i+1}|,\cdots,|x_{n-1}|) < \mathrm{ulp}(y_{i-1})$.
\end{proof}		

We now analyze several conditions under which the VecSum output is F-nonoverlapping. These conditions are similar to those in \cite{Popescu2017Phd}, but yield stronger results.

\begin{theorem}
	If $\{x_i\}_{i=0}^{n-1}$ overlaps by at most $d\leq p-2$ bits (i.e., the exponent of $x_i$ is larger than that of $x_{i+1}$ by at least $2$), we do not require three consecutive terms to have no overlapping digits, and ties are broken by the ties-to-even rule, then $\{y_i\}=\mathrm{VecSum}(\{x_i\})$ is F-nonoverlapping. Moreover, Fast2Sum can be used in VecSum, and there is no $0 \leq j < i <n$ such that $ |y_{j}|=2^{k+1}$ and $|y_{i}|=2^{k}$ with $k\in\mathbb{Z}$.
\end{theorem}
\begin{proof}
	Suppose $x_{i-1}=\frac M {2^{p-1}} \cdot 2^k \geq 2^k$ with $2^{p-1}\leq M < 2^p$ and $M\in\mathbb{Z}$, then $|x_j| \leq \frac {2^p-1} {2^{p-1}} 2^{k-2(j-(i-1))} < 2^{k-1-2(j-i)}$. $|t_i|=|\mathrm{RORS}(x_i,x_{i+1},\cdots,x_{n-1})|\leq |x_{i-1}| \leq 2^{k-1}\mathrm{RORS}(1,2^{-2},\cdots,2^{-2(n-1-i)} )\leq \frac 4 3 2^{k-1}$, which implies $|x_{i-1}| > |t_i| $ and $\mathrm{ulp}(x_{i-1})\geq 2\mathrm{ulp}(t_{i})$. By Lemma \ref{lem:Fbasecond}, $\{y_i\}$ is F-nonoverlapping.
	
	Suppose $|y_{i}|=2^{k}$ with $k\in\mathbb{Z}$. If $|y_i|=\frac 1 2 \mathrm{ulp}(t_{i-1})$, then, by the ties-to-even rule, the last significant digit of $|t_{i-1}|$ is zero, so that $\mathrm{uls}(t_{i-1})\geq 2\mathrm{ulp}(t_{i-1})=4|y_i|$; otherwise $|y_i|\leq \frac 1 4 \mathrm{ulp}(t_{i-1})$. In summary, $\mathrm{uls}(t_{i-1})\geq 4|y_i|$. Also, $\mathrm{uls}(x_k)\geq \mathrm{ulp}(x_{k})\geq 2\mathrm{ulp}(t_{i-1})\geq 4|y_i|$ for $0\leq k <i-1$. Therefore $|y_j|\geq \mathrm{uls}(y_j)\geq 4|y_i|$ or $y_j=0$ for $0\leq j <i$.
	
\end{proof}

\begin{theorem}
	If $\{x_i\}_{i=0}^{m-1}$ overlaps by at most $d\leq p-1$ bits (i.e., the exponent of $x_i$ is larger than that of $x_{i+1}$ by at least $1$), and we do not require three consecutive terms to have no overlapping digits, then $\{y_i\}=\mathrm{VecSum}(\{x_i\})$ is F-nonoverlapping, and Fast2Sum can be used in VecSum.
\end{theorem}
\begin{proof}
	Suppose $|x_0|=M_0\cdot 2^k$ with $2^{p-1}\leq M_0 < 2^p$, and $|x_i|=M_i\cdot 2^{k_i}$ with $2^{p-1}\leq M_i < 2^p$, $k_i\leq k-i$, and $k_{i+1}<k_i$. It suffices to prove that the exponent of $\mathrm{RORS}(x_1,x_2,\cdots, x_n)$ does not exceed that of $x_0$, i.e.,  $|\mathrm{RORS}(x_1,x_2,\cdots, x_n)| \leq (2^p-1)\cdot 2^k$.
	\begin{align*}
	|\mathrm{RORS}(x_1,x_2,\cdots, x_n)| &= |\mathrm{RN}(x_1 + \mathrm{RN}(x_2 + \cdots +\mathrm{RN}(x_{n-2}+x_{n-1})\cdots)|  \\
	&\leq |\mathrm{RN}(|x_1| + \mathrm{RN}(|x_2| + \cdots +\mathrm{RN}(|x_{n-2}|+|x_{n-1}|)\cdots)| \\
	&\leq \mathrm{RN}((2^p-1)\cdot 2^{k-1} + \mathrm{RN}((2^p-1)\cdot 2^{k-2}  + \cdots\\ &\quad +\mathrm{RN}((2^p-1)\cdot 2^{k-(n-2)} +(2^p-1)\cdot 2^{k-(n-1)} )\cdots))\\
	\end{align*}
	We will prove $\mathrm{RN}((2^p-1)\cdot 2^{k-1} + \mathrm{RN}((2^p-1)\cdot 2^{k-2}  + \cdots +\mathrm{RN}((2^p-1)\cdot 2^{k-(n-2)} +(2^p-1)\cdot 2^{k-(n-1)} )\cdots)) \leq  (2^p-1)\cdot 2^k$ by induction.
	When $n=2$, clearly $|\mathrm{RORS}(x_1,x_2,\cdots, x_n)|= |x_1| \leq (2^p-1)\cdot 2^{k-1} < (2^p-1)\cdot 2^{k}$.
	If $n > 2$, the induction hypothesis gives $\mathrm{RORS}(x_2,\cdots, x_n)\leq (2^p-1)\cdot 2^{k-1}$, hence $\mathrm{RORS}(x_1,x_2,\cdots, x_n)\leq \mathrm{RN}( (2^p-1)\cdot 2^{k-1} +(2^p-1)\cdot 2^{k-1})=(2^p-1)\cdot 2^{k}.$
	
	This implies that the exponent of $t_i$ does not exceed that of $x_{i-1}$; by Lemma \ref{lem:Fbasecond}, $\{y_i\}$ is F-nonoverlapping.
	
	Under this condition, one cannot guarantee that there is no $0 \leq  j < i <n$ such that $ |y_j|=2^{k+1}$ and $|y_{i}|=2^{k}$ with $k\in\mathbb{Z}$.
	For example, $\{x_i\}=(-2,1-2^{1-p},2^{-1}+2^{-2}, 2^{-2}+2^{-p})$ yields the output $\{y_i\}=(-2^{1-p},0,0,2^{-p})$.	
	
\end{proof}

\begin{algorithm}[H]\caption{MultiVecSum$(x_0,x_1,\cdots,x_{m-1})$}\label{algo:multivecsum}
	\begin{algorithmic}[1]
		\REQUIRE $\{x_i\}_{i=0}^{m-1}$ is F-nonoverlapping; the first two nonzero terms $x_i,x_j$ with $i<j$ do not satisfy $|x_i|=2^{k+1}$, $|x_j|=2^{k}$, $x_ix_j<0$ for $k\in\mathbb{Z}$.
		\ENSURE $\{y_i\}_{i=0}^{n-1}$ is ulp-nonoverlapping and has no interleaving zeros.
		\STATE $\{s_i^{(0)}\}_{i=0}^{m-1} \leftarrow \{x_i\}_{i=0}^{m-1}$
		\FOR {$j\leftarrow 0 $ \TO $n-2$} 
		\STATE $\{s_i^{(j+1)}\}_{i=j}^{m-1}\leftarrow\mathrm{VecSum}(\{s_i^{(j)}\}_{i=j}^{m-1})$
		\STATE $y_j \leftarrow s_j^{(j+1)}$
		\ENDFOR
		\STATE $y_{n-1}\leftarrow \mathrm{RORS}(\{s_i^{(n-1)}\}_{i=n-1}^{m-1})$
		\RETURN  $y_0,y_1,\cdots,y_{n-1}$
	\end{algorithmic}
\end{algorithm}

This completes the proof of the correctness of Algorithm \ref{algo:multivecsum}.

\begin{theorem}
	Let $\{z_i\}_{i=0}^{n-1}$ be the result of merging the expansions $\{x_i\}$ and $\{y_i\}$. According to \cite{Popescu2017Phd}, $\{f_i\}_{i=0}^{n-1}\leftarrow\mathrm{VecSum}(\{z_i\})$ is strictly $\mathcal{S}$-nonoverlapping, or $f_i=f_{i+1}=\pm 2^k$ for some $k\in\mathbb{Z}$.
	We show that $\{g_i\}_{i=0}^{n-1} \leftarrow \mathrm{VecSum}(\{f_i\})$ is F-nonoverlapping, that Fast2Sum can be used in VecSum, and that there is no $0 <  i <n$ such that $ |f_{i-1}|=2^{k+1}$ and $|f_{i}|=2^{k}$ with $k\in\mathbb{Z}$.
\end{theorem}
\begin{proof}
	First, if $|e_i|\leq 2^k$, then $|t_i|=|\mathrm{RORS}(e_i,e_{i+1},\cdots, e_{n-1})|\leq |\mathrm{RORS}(2^k,2^k,2^{k-1},2^{k-1},2^{k-2},2^{k-2},\cdots)| \leq 2^{k+2}$ (by induction).
	$(t_{i-1},f_i) \leftarrow \mathrm{Fast2Sum}(e_{i-1},t_i)$,  	$(t_{i-2},f_{i-1}) \leftarrow \mathrm{Fast2Sum}(e_{i-2},t_{i-1})$.

	\begin{itemize}
		\item  If $|f_{i-2}|=2^k$, $|f_{i-1}|\leq 2^k$, $|f_{i}|\leq 2^{k-1}$, $|t_i|=|\mathrm{RORS}(f_i,f_{i+1},\cdots,f_{n-1})|\leq 2^{k+1}$,  $|t_{i-1}|=|\mathrm{RN}(f_{i-1}+t_{i-1})|\leq 3\cdot 2^k$.
		$|g_i|\leq \frac 1 2 \mathrm{ulp}(t_{i-1})\leq \frac 1 2 \mathrm{uls}(t_{i-1})$, and $\mathrm{uls}(f_{i-2})=2^k \geq \mathrm{ulp}(t_{i-1})$, so $|g_i|\leq \frac 1 2 \mathrm{uls}(f_{i-2})$; therefore $|g_i|\leq \frac 1 2 \mathrm{uls}(g_{i-1})$.
		Since $\mathrm{uls}(f_{i-2})=2^k \geq \mathrm{ulp}(t_{i-1})$, Fast2Sum is exact by Lemma \ref{lem:fast2sumextra}.
		
		$\mathrm{uls}(f_{i-2})=2^k > 4\cdot(2^{k+1}\cdot 2^{-p} )\geq4\cdot \frac {1} {2} \mathrm{ulp}(t_{i-1})$.
		
		\item  $|t_{i-1}|\leq 2\mathrm{uls}(f_{i-2})$. If $|f_{i-2}|\neq2^k$, then $|f_{i-2}|\geq 3\mathrm{uls}(f_{i-2})$; since $ |t_{i-1}|<|f_{i-2}|$, $\mathrm{Fast2Sum}(f_{i-2},t_{i-1})$ is exact, and $|g_i|\leq \frac 1 2 \mathrm{uls}(g_{i-1})$.
		
		If $|f_{i-2}|=(1+2)\mathrm{uls}(f_{i-2})$, then $\mathrm{uls}(f_{i-2})\geq 2^{p-2}\mathrm{ulp}(f_{i-2})\geq 2^{p-2}\mathrm{ulp}(t_{i-1})\geq 4\cdot \frac 1 2 \mathrm{ulp}(t_{i-1})$.
		
		Otherwise $|f_{i-2}|\geq (1+4)\mathrm{uls}(f_{i-2})  > 2|t_{i-1}|$, and $\mathrm{ulp}(f_{i-2}) \geq 2 \mathrm{ulp}(t_{i-1}) = 4\cdot \frac 1 2 \mathrm{ulp}(t_{i-1})$. 
	\end{itemize}
	By the ties-to-even rule, $\mathrm{uls}(t_{i-1})\geq 4 \cdot \frac 1 2\mathrm{ulp}(t_{i-1})$.  Therefore, for all $j<i$, we always have $\mathrm{uls}(g_j)\geq 4\cdot \frac  1 2\mathrm{ulp}(t_{i-1})$, or $g_j=0$. 
	
	The output $f_0$ may be zero or inaccurate, in which case at least one more $\mathrm{VecSum}(\{f_i\}_{i=1}^{n-1})$ is required, e.g., $\{x_i\}=(1,2^{-10}-1, 2^{-p-2},2^{-p-2})$, $\{f_i\}=(2^{-10},0,2^{-p-1},0)$.	
\end{proof}

\begin{theorem}
	For $p\geq9$, suppose $\{x_i\}$ and $\{y_i\}$ are ulp-nonoverlapping and either $x_0y_0>0$, or $x_0y_0<0$ and $\frac {\max(|x_0|,|y_0|)} {\min(|x_0|,|y_0|)} \geq 2(1+10\mathbf{u})$. If $\{f_i\}\leftarrow \mathrm{VecSum}(\mathrm{Merge}(\{x_i\},\{y_i\}))$, then $f_0$ is accurate, i.e., $|\mathrm{RORS}(f_1,f_2,\cdots)|<\mathrm{ulp}(f_0)$.
\end{theorem}
\begin{proof}
	Let $ \{z_i\} \leftarrow \mathrm{Merge}(\{x_i\},\{y_i\})$. 
	Without loss of generality, assume $|x_0|\geq |y_0|$.
	Let $\eta = \mathrm{ulp}(x_0)$ and $v = 2^{1-p}$. 
	\begin{itemize}
		\item $|y_0|\leq |x_1| \leq \mathrm{ulp}(x_0)$.  $ z_0 = x_0$, $ (|z_1|,|z_2|,|z_3|,|z_4|,\cdots,|z_{2n-3}|,|z_{2n-2}|,|z_{2n-1}|) \leq (\eta,\eta, v\eta, v\eta,\cdots,v^{n-2}\eta,v^{n-2}\eta,v^{n-1}\eta )$ in the element-wise sense.
		\begin{align*}
		|t_{2n-2}| &\leq \mathrm{RN}(v^{n-2}\eta+v^{n-1}\eta)=(1+v)v^{n-2}\eta \\
		|t_{2n-3}| &\leq \mathrm{RN}((1+v)v^{n-2}\eta+v^{n-2}\eta)\leq 2(1+v)v^{n-2}\eta \\
		|t_{2n-4}| &\leq \mathrm{RN}( 2(1+v)v^{n-2}\eta+v^{n-3}\eta)=(1+2v)v^{n-3}\eta \\
		|t_{2n-5}| &\leq \mathrm{RN}((1+2v)v^{n-3}+v^{n-3}\eta)=2(1+v)v^{n-3}\eta\\
		\cdots
		\end{align*}
		By induction we have $|t_{2n-2k}|\leq (1+2v)v^{n-1-k}\eta$, $|t_{2n-(2k+1)}|\leq 2(1+v)v^{n-1-k}\eta$, $1\leq k \leq n-1$.
		
		$|t_1|\leq 2(1+v)\eta$, $|t_2|\leq (1+2v)\eta$.
		$|f_2|\leq \frac 1 2 \mathrm{ulp}(t_1)\leq \frac v 2 |t_1|\leq (1+v)v\eta$.  In the previous proof, we have shown $|\mathrm{RORS}(f_3,f_4,\cdots) |\leq 3|f_2|$, hence $|\mathrm{RORS}(f_2,f_3,\cdots)| \leq 4|f_2| \leq 4(1+v)v\eta$.
		$|f_1|\leq \frac 1 2 \mathrm{ulp}(f_0)$, 
		$ |\mathrm{RORS}(f_2,f_3,\cdots)| \leq 4(1+v)v\eta \leq 8(1+v)v\mathrm{ulp}(f_0)$. $|\mathrm{RORS}(f_1,f_2,\cdots )| <  \mathrm{ulp}(f_0)$ for $p\geq 6$.
		\item $x_0y_0>0$, $|y_0| >  |x_1|$, $z_0=x_0$, $z_1=y_0$. We still have $( |z_2|,|z_3|,|z_4|,\cdots,|z_{2n-3}|,|z_{2n-2}|,|z_{2n-1}|) \leq (\eta,\eta, v\eta, v\eta,\cdots,v^{n-2}\eta,v^{n-2}\eta )$ in the element-wise sense.
		\begin{align*}
		|t_{2n-2}|&\leq\mathrm{RN}(v^{n-2}\eta+v^{n-2}\eta)=2v^{n-2}\eta\\
		|t_{2n-3}|&\leq\mathrm{RN}(2v^{n-2}\eta+v^{n-3}\eta)=(1+2v)v^{n-3}\eta\\
		|t_{2n-4}|&\leq\mathrm{RN}((1+2v)v^{n-3}\eta+v^{n-3}\eta)=2(1+v)v^{n-3}\eta\\
		|t_{2n-5}|&\leq\mathrm{RN}(2(1+v)v^{n-3}\eta+v^{n-4}\eta)=(1+2v)v^{n-4}\eta\\
		\end{align*}
		By induction, $ |t_{2n-2k}|\leq 2(1+v)v^{n-1-k}\eta$, $|t_{2n-(2k+1)}|\leq (1+2v)v^{n-2-k}\eta$.
		$|t_2|\leq 2(1+v)\eta$. $|f_3|\leq \frac 1 2 \mathrm{ulp}(t_2)\leq \frac v 2 |t_2|$, $|\mathrm{RORS}(f_3,f_4,\cdots) |\leq 4|f_3| \leq 2\mathrm{ulp}(t_2)$.
		\begin{itemize}
			\item $y_0 t_1<0$, i.e., $x_0t_1<0$. This implies $|y_0|<|t_2|$, $|t_1|\leq |t_2|$,  $|f_2|\leq \frac 1 2 \mathrm{ulp}(t_1) \leq \frac 1 2 \mathrm{ulp}(t_2)$. $|\mathrm{RORS}(f_2,f_3,\cdots)|\leq \frac 5 2 \mathrm{ulp}(t_2)\leq \frac 5 2 \cdot v \cdot 2(1+v)\eta=5v(1+v)\eta$. Since $|t_1|\leq |t_2|\leq 2(1+v)\mathrm{ulp}(x_0)$, $\mathrm{ulp}(f_0)=\mathrm{ulp}(\mathrm{RN}(x_0+t_1))\geq \frac 1 2 v$. $|\mathrm{RORS}(f_1,f_2,\cdots )| <  \mathrm{ulp}(f_0)$ for $p\geq 5$.
			\item $y_0t_1\geq 0$, i.e., $x_0t_1\geq 0$.  $e_{y_0}\leq e_{x_0}$, $e_{y_0}-1\leq e_{t_1}\leq e_{y_0}+1$. $e_{f_0}\geq e_{x_0}$.  $|\mathrm{RORS}(f_3,f_4,\cdots) |\leq 4|f_3| \leq 2\mathrm{ulp}(t_2)\leq 4v(1+v)\eta \leq \frac 3 {16} \eta \leq \frac 3 {16} \mathrm{ulp}(f_0) $ for $p\geq 6$.
			\begin{itemize}
				\item $e_{t_1}=e_{y_0}-1\leq e_{f_0}-1$, $|f_2|\leq \frac 1 2 \mathrm{ulp}(t_1)\leq \frac 1 4 \mathrm{ulp}(f_0)$, $|\mathrm{RORS}(f_1,f_2,\cdots)|< \mathrm{ulp}(f_0)$.
				\item $e_{t_1}=e_{y_0}$. If $e_{y_0}=e_{x_0}$, then $e_{f_0}=e_{x_0}+1=e_{t_1}+1$. The argument is the same as above.
				Otherwise, if $e_{y_0}<e_{x_0}$, the same argument applies. 
				\item $e_{t_1}=e_{y_0}+1$. If $e_{y_0} <e_{x_0}$, $e_{t_1}<e_{f_0}$. Otherwise if $e_{y_0}=e_{x_0}, e_{t_1}=e_{y_0}+1$. $|t_1|=|\mathrm{RN}(y_0+t_2)$. $|t_2|\leq 2(1+v)\eta$. Suppose $|x_0|=M_x 2^k, |y_0|=M_y 2^k, 2^{p-1}\leq M_x,M_y<2^p$.  Because  $e_{t_1}=e_{y_0}+1$, $M_x\geq M_y \geq 2^p-2$.
				\begin{itemize}
					\item $M_x\geq M_y=2^p-2$, then $|t_2|\geq \frac 3 2 \eta$, $|t_1|=2^p2^k$. If $M_x=2^p-1$, then  $|f_0|\geq\mathrm{RN}( (2^p-1)2^k+2^p2^k )=2^{p+1}2^k$, $e_{f_0}=e_{t_1}+1$. 
					Otherwise $M_x=2^p-2$, $x_0+t_1=f_0$, i.e., the addition is exact, $f_1=0$. $|\mathrm{RORS}(f_1,f_2,f_3,\cdots)|\leq \frac 1 2 \mathrm{ulp}(t_1)+|\mathrm{RORS}(f_3,f_4,\cdots)|\leq \frac 1 2 \mathrm{ulp}(f_0)+\frac 3 {16}\mathrm{ulp}(f_0)<\mathrm{ulp}(f_0)$.
					\item $M_x=M_y=2^p-1$, $2^{p+1}2^k>|t_1|\geq 2^p2^k$, $|f_0|\geq\mathrm{RN}( (2^p-1)2^k+2^p2^k )=2^{p+1}2^k$, $e_{f_0}=e_{t_1}+1$. 
				\end{itemize}
			\end{itemize}
		\end{itemize}
		
		\item $x_0y_0<0$, $|y_0| >  |x_1|$, $z_0=x_0$, $z_1=y_0$. $|t_0|\leq 2(1+v)\eta \leq 2v(1+v)|x_0|$, $|t_1|\leq (|y_0| + 2v(1+v)|x_0|)(1+\mathbf{u})$.
		First we show that it suffices to prove $|x_0|\geq 2|t_1|$.  In this case $ |f_0| \geq |t_1|$.
		\begin{itemize}
				\item $e_{f_0} > e_{t_1}$, $|f_1|\leq \frac 1 2 \mathrm{ulp}(f_0)$, $|f_2|\leq \frac 1 4 \mathrm{ulp}(f_0)$, and $\{f_i\}$ is strictly F-nonoverlapping, so $|\mathrm{RORS}(f_1,f_2,\cdots)|< \mathrm{ulp}(f_0)$.
			\item $e_{f_0} = e_{t_1}$. $\mathrm{ulp}(x_0)=2^ke_{\mathrm{ulp}(t_1)}=2^ke_{\mathrm{ulp}(f_0)}$ with $k\in\mathbb{N}$; thus $f_1$ must be a multiple of $\mathrm{ulp}(f_0)$ or zero, and $|f_1|\leq \frac 1 2 \mathrm{ulp}(f_0)$. Hence $f_1$ can only be zero. Also $|f_2|\leq \frac 1 2\mathrm{ulp}(f_0)$, and $\{f_i\}$ is strictly F-nonoverlapping, so $|\mathrm{RORS}(f_1,f_2)| < \mathrm{ulp}(f_1)$.
		\end{itemize}
		$|x_0|\geq 2|t_1| \Leftrightarrow |x_0| \geq 2(|y_0| + 2v(1+v)|x_0|)(1+\mathbf{u})=2(1+\mathbf{u})(|y_0|+4\mathbf{u}(1+2\mathbf{u})|x_0|)$, $|x_0| \geq 2 \frac {1+\mathbf{u}}  {1-8\mathbf{u}(1+\mathbf{u})(1+2\mathbf{u})}|y_0|$. 
		$ (1-8\mathbf{u}(1+\mathbf{u})(1+2\mathbf{u}))\cdot(1+8\frac 1 4\mathbf{u})= 1+\mathbf{u}(\frac 1 4-90\mathbf{u}-241\mathbf{u}^2-132\mathbf{u}^3) >1 $ when $p\geq 9, \mathbf{u}\leq 2^{-9}$.
		$|x_0|\geq 2(1+\mathbf{u})(1+8\frac 1 4 \mathbf{u})|y_0| = 2(1+9\frac 1 4 \mathbf{u} + 8\frac 1 4 \mathbf{u}^2)|y_0|$.	 $|x_0|\geq 2(1+10\mathbf{u})|y_0|$.	
		Counterexample: $x_0=(2^p-1)2^k$, $y_0=-(2^p-1)2^{k-1}$, $\eta = 2^k$, $t_2=-2^k$. The output is $(f_0,f_1,f_2)=((2^p-2)2^{k-1},0,-2^{k-1})$ with $|f_2|=\mathrm{ulp}(f_0)$.
	\end{itemize}
	
\end{proof}

\begin{remark}
	The VecSum of a merged expansion may not be F-nonoverlapping. For example, for the input
	\{0x1.FFFFFFFFFFFFFp53,  0x1.FFFFFFFFFFFFFp53,	0x1.FFFFFFFFFFFFFp0,   0x1.FFFFFFFFFFFFFp0,	0x1.FFFFFFFFFFFFFp-53, 0x1.FFFFFFFFFFFFFp-53\}, where $\{x_i\}$ and $\{y_i\}$ are $\mathcal{P}$-nonoverlapping, the output is
	\{0X1P+55 -0X1P+1 +0X1P+1 -0X1P-52 +0X1.FFFFFFFFFFFFEP-53 +0X0P+0\}. 
	Moreover, even though $e_0=0X1P+55$ is accurate enough, $\mathrm{VecSum}(\{f_i\}_{i=1}^{n-1})=(0,0,-0x1p-104,0,0)$, i.e., another $\mathrm{VecSum}(\{f_i\}_{i=1}^{n-1})$ is required.
	If we prepend $\{-0x1p106-0x1p55, 0x1p106\}$, i.e., use the input \{-0x1p106-0x1p55, 0x1p106,
	0x1.FFFFFFFFFFFFFp53 ,  0x1.FFFFFFFFFFFFFp53,
	0x1.FFFFFFFFFFFFFp0,   0x1.FFFFFFFFFFFFFp0,
	0x1.FFFFFFFFFFFFFp-53, 0x1.FFFFFFFFFFFFFp-53\}, the result is more interesting. 
	
	After one VecSum, the output is: \{0x0P+0, +0X0P+0, +0X0P+0, -0X1P+1, +0X1P+1, -0X1P-52, +0X1.FFFFFFFFFFFFEP-53, +0X0P+0\};\\
	after two VecSums, the output is: \{0X0P+0, +0X0P+0, +0X0P+0, +0X0P+0, +0X0P+0, -0X1P-104, +0X0P+0, +0X0P+0\};\\
	after three VecSums, the output is: \{-0X1P-104, +0X0P+0, +0X0P+0, +0X0P+0, +0X0P+0, +0X0P+0, +0X0P+0, +0X0P+0\}.\\
	
	The result can be $\{0,\cdots,0, 2^k, -2^k,\cdots\}$, and at least three VecSums are required to ensure that the first term is exact.
	However, such cases are rare. From the proof in \cite{Popescu2017Phd} we can see that this happens only when $z_i$ and $z_{i+1}$ overlap in $p$ digits, $\mathrm{ulp}(t_{i+1})=2\mathrm{ulp}(z_i)$, and $|f_{i+1}|=|f_{i+2}|=\mathrm{ulp}(z_i)$. Because of the ulp-nonoverlapping condition, we know that $|t_{i+1}|\leq 2(1+\epsilon) \mathrm{ulp}(z_{i+1})$. To obtain $f_{i+2}=\mathrm{ulp}(z_i)$, we must have $|z_i|=|z_{i+1}|=(2^p-1)2^k$ and $|t_{i+1}|\geq 2\mathrm{ulp}(z_i)$. But in this case, if $z_iz_{i+1}<0$, then $f_{i+1}f_{i+2}>0$, whereas $f_{i+1}f_{i+2}<0$ can occur only when $z_iz_{i+1}>0$. To cancel $|t_i|=2^p2^{k+1}$, at least two nearby floating-point numbers with opposite signs must precede.
	
	In summary, if there is no cancellation in the input, $f_0$ is already accurate after the first step following the merge, but $\{f_i\}_{i=1}^{n-1}$ still needs another VecSum to be converted into a form acceptable to MultiVecSum. In general, two VecSums are required before MultiVecSum.
\end{remark}

\section{Review of the QDAdd Algorithm}

\begin{figure}[H]
	\centering
	\includegraphics[width=0.7\linewidth]{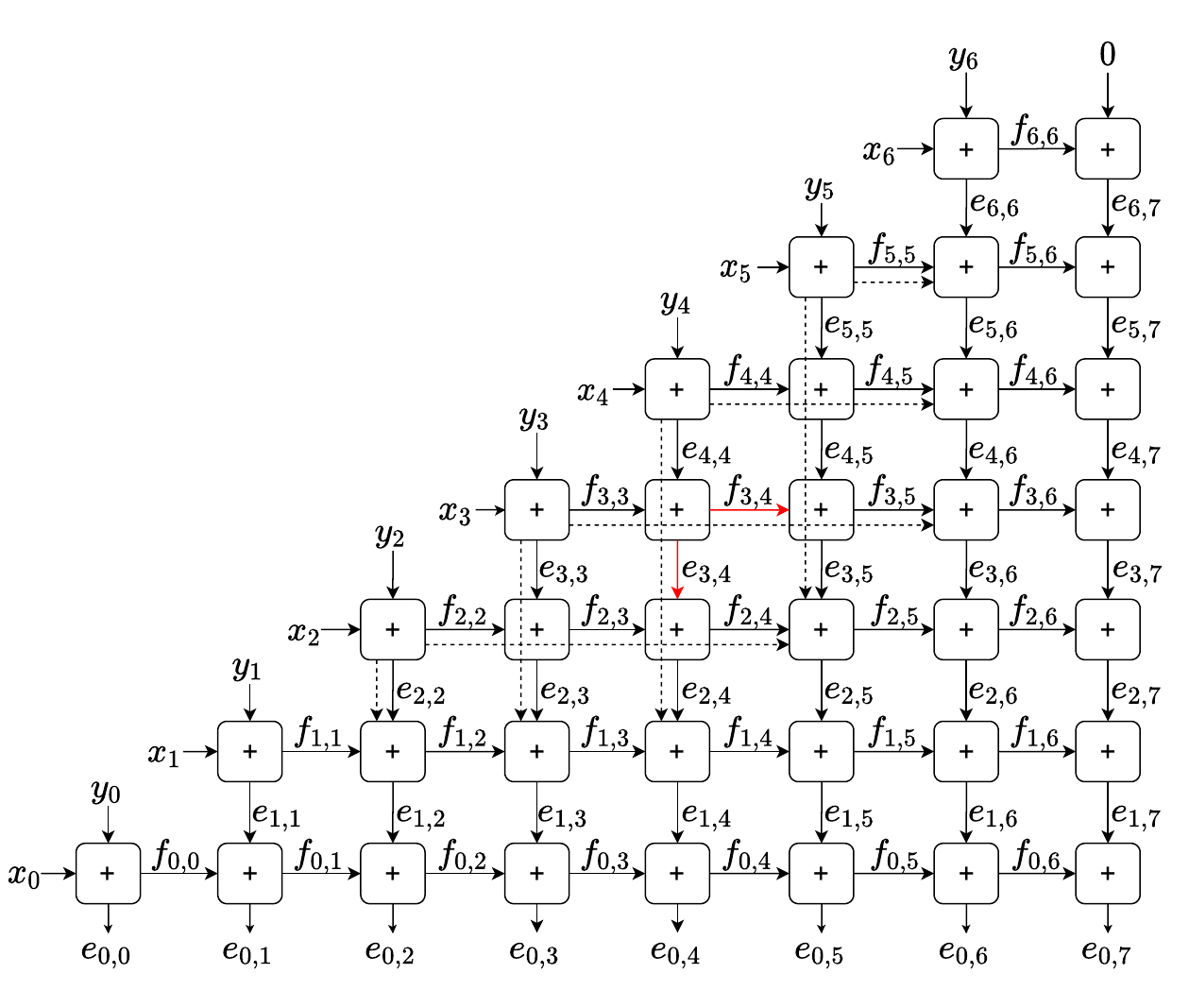}
	\caption{Bailey's addition algorithm}
	\label{fig:baileyadd}
\end{figure}

The assertion in \cite{Popescu2017Phd} that ``All the terms are of the same order, and it can be shown by induction that $e_{r-1}$ ($f_{6,6}$ in Figure \ref{fig:baileyadd}) is the largest term in this addition'' is misleading. As a counterexample, let $|y_0|<\frac 1 2 \mathrm{ulp}(x_1)$ and $|y_1| < \frac 1 2 \mathrm{ulp}(y_0)$; then $e_{1,1}=x_1$, $f_{1,1}=y_1$, $e_{0,0}=x_0$, $f_{0,0}=y_0$, $f_{0,1}=y_0$, and $ |f_{0,1}|=|y_0|>|y_1|=|f_{1,1}|$.

Therefore, the error estimate given in that paper is incorrect; we give a correct (though not necessarily tight) error estimate here.
\begin{theorem}
	Consider the BaileyAdd diagram in Figure \ref{fig:baileyadd}. When the inputs are ulp-nonoverlapping, the absolute error of the Bailey Add algorithm before renormalization is at most $\gamma_{n-1} \mathbf{u}^n(|x_0|+|y_0|)(3^n(1+\gamma_n)-2^n)$.
\end{theorem}
\begin{proof}
	Every $e_{i,j}$ with $i\neq j$ originates from some $e_{k,k}$ and $f_{k,k}$, and $|f_{k,k}|\leq \mathbf{u}|e_{k,k}|$; hence the upper bound on $|e_{i,j}|$ is determined by the $e_{k,k}$ that can propagate to $e_{i,j}$. This propagation is one-directional: it can only go right or downward. Let the horizontal distance be $\Delta h$ and the vertical distance be $\Delta v$. Then $e_{k,k}$ must undergo $\Delta h+\Delta v$ floating-point additions to reach $e_{i,j}$; taking floating-point errors into account, this introduces a scaling factor $(1+\mathbf{u})^{\Delta h + \Delta v}$, and since the operations along the horizontal direction extract the floating-point error, an additional factor $\mathbf{u}^{\Delta h}$ is required. The total path length is $\Delta h+\Delta v$, and each path amounts to choosing $\Delta h$ rightward steps (or $\Delta v$ downward steps) out of $\Delta h+\Delta v$ steps; hence there are $C_{\Delta h + \Delta v}^{\Delta v}$ such paths, where $\Delta h = j-k$ and $\Delta v = k-i$. Therefore
	\begin{align*}
	|e_{i,j}| \leq \Sigma_{k=i}^{j} C_{j-i}^{k-i}|e_{k,k}| \mathbf{u}^{j-k}(1+\mathbf{u})^{j-i}
	\end{align*}
	
	Assume $(2\mathbf{u})^i 2^{-d_i} |x_0| \leq |x_i| \leq (2\mathbf{u})^i |x_0|$ and $(2\mathbf{u})^i 2^{-d_i} |y_0| \leq |y_i| \leq (2\mathbf{u})^i |y_0|$. Then $|e_{k,k}|  \leq |x_k+y_k|(1+\mathbf{u}) \leq (2\mathbf{u})^k(|x_0|+|y_0|)(1+\mathbf{u})$.
	
	Before renormalization, the error of the Bailey addition comes entirely from the accumulation in the rightmost layer, namely $\Sigma_{i=0}^{n-1} f_{i,n-1}$, where $|f_{i,n-1}|\leq \mathbf{u} |e_{i,n-1}|$. Hence the accumulated sum satisfies
	\begin{align*}
	|S_1 | & = \Sigma_{i=0}^{n-1}|f_{i,n-1}| \\
	& \leq \mathbf{u} \Sigma_{i=0}^{n-1} |e_{i,n-1}| \\
	& \leq \Sigma_{i=0}^{n-1}(1+\mathbf{u})^{n-1-i} \left (\Sigma_{k=i}^{n-1} C_{n-1-i}^{k-i}|e_{k,k}|\mathbf{u}^{n-1-k}\right )\\
	&\leq \mathbf{u}(|x_0|+|y_0|)\Sigma_{i=0}^{n-1}(1+\mathbf{u})^{n-i} \left(  
	\Sigma_{k=i}^{n-1} C_{n-1-i}^{k-i} 2^k\mathbf{u}^{n-1}
	\right )\\
	&= \mathbf{u}^n(|x_0|+|y_0|)\Sigma_{i=0}^{n-1}(1+\mathbf{u})^{n-i} \left(  		            \Sigma_{k=i}^{n-1} C_{n-1-i}^{k-i} 2^k
	\right )\\
	&\leq \mathbf{u}^n(|x_0|+|y_0|) ((3+3\mathbf{u})^n-2^n) \\
	&= \mathbf{u}^n(|x_0|+|y_0|) (3^n(1+\mathbf{u})^n-2^n)\\
	&\leq \mathbf{u}^n(|x_0|+|y_0|)(3^n+1-2^n)
	\end{align*}
	
	In the last step, $3^n(1+\mathbf{u})^n\leq 3^n+3^n\gamma_n = 3^n+\frac {3^n n 2^{-p}} {1-n2^{-p}}$; for the commonly used ranges $n\leq p/2$ and $p\geq 24$, we always have $3^n\gamma_n \leq 1$.
	
	The absolute error before renormalization is less than $\gamma_{n-1}|S_1|\leq \gamma_{n-1} \mathbf{u}^n(|x_0|+|y_0|)(3^n+1-2^n)$. When $n\mathbf{u}\ll 1$, this bound is approximately $n(3^n-2^n+1)\mathbf{u}^{n+1}(|x_0|+|y_0|)$, which remains roughly of the order of $O(\mathbf{u}^{n+1}(|x_0|+|y_0|))$ for moderate $n$.
\end{proof}

\begin{theorem}\label{thm:baileyulsulp}
	Assume $n\geq 3$, $p\geq 8$, and $2^p>2^{n-1}(n+1)+6$. Then, in the BaileyAdd diagram of Figure \ref{fig:baileyadd}, $\mathrm{uls}(e_{0,j})\geq \mathrm{ulp}(\mathrm{RORS}(e_{0,j+1},e_{0,j+2},\cdots,e_{0,n} ))$.	
\end{theorem}
Without the above theorem, one would regard $\{e_{0,j}\}_{j=0}^n$ as random and would need repeated applications of VecSum to obtain an F-nonoverlapping result. With the above theorem, a single $\mathrm{VecSum}(\{e_{0,j}\}_{j=0}^n)$ suffices. The counterexample given in \cite{Popescu2017Phd} merely violates $\mathrm{ulp}(e_{0,j})\geq \mathrm{ulp}(\mathrm{RORS}(e_{0,j+1},e_{0,j+2},\cdots,e_{0,n} ))$, while $\mathrm{uls}(e_{0,j})\geq \mathrm{ulp}(\mathrm{RORS}(e_{0,j+1},e_{0,j+2},\cdots,e_{0,n} ))$ still holds; it is the latter that is the essential condition.
In the original source of the Bailey addition algorithm \cite{Bailey2000QD}, the future-work section already pointed out that the input before renormalization does not satisfy the conditions of the then-known VecSum+VecSumErrBranch algorithm. \cite{Popescu2017Phd} solved this problem by applying VecSum+VecSumErrBranch repeatedly, a renormalization algorithm suitable for random inputs, but possibly because of performance drawbacks, later returned to the original method of the QD library and called it the quick-and-dirty algorithm; the conclusion here shows that this method is not dirty at all. 

To prove the above theorem, we first establish several lemmas.

\begin{lemma}\label{lem:2sumcompulp}
	Let $a,b\in\mathbb{F}_p$, $s=\mathrm{RN}(a+b)$, and $s+t=a+b$. Then $|t|\leq|a|$, $|t|\leq |b|$, and $\mathrm{ulp}(t)\leq \min(\mathrm{ulp}(a),\mathrm{ulp}(b))$.
\end{lemma}
\begin{proof}
	Let $r=a+b$. By the definition of RN, $s$ is the floating-point number nearest to the real number $r$, and $|t|$ is the distance between $s$ and $r$. Since $|a|$ and $|b|$ are the distances from $b$ and $a$ to $r$, respectively, we have $|t|\leq|a|$ and $|t|\leq |b|$.
\end{proof}

\begin{lemma}
	Let $x$ and $y$ be two floating-point numbers, and let $z = \max(|x|,|y|)$, $w = \min(|x|,|y|)$, and $s=\mathrm{RN}(x+y)$. Then $\mathrm{uls}(s)\geq \frac 1 2\mathrm{ulp}(z)$ or $s=0$.
\end{lemma}
\begin{proof}
	If $r(a,b)<\frac 1 2 $, then $|s|\geq \frac z 2$ and $\mathrm{uls}(s)\geq \mathrm{ulp}(s)\geq \frac 1 2 \mathrm{ulp}(z)$; otherwise, if $r(a,b)\geq \frac 1 2 $, then $|w|\geq \frac 1 2|z|$ and $\mathrm{ulp}(w)\geq \frac 1 2 \mathrm{ulp}(z)$, so $\mathrm{uls}(s)\geq \min(\mathrm{uls}(z),\mathrm{uls}(w))\geq\min(\mathrm{ulp}(z),\mathrm{ulp}(w))\geq\frac 1 2 \mathrm{ulp}(z) $. 
\end{proof}
\begin{lemma}\label{lem:binsel}
	Let $a,b\in \mathbb{F}$, $\mathrm{RN}(a+b)=s$. Then at least one of $|a|$ and $|b|$ is greater than or equal to $|s|/2$.
\end{lemma}
\begin{proof}
	The case $s=0$ is trivial. For $s\neq 0$, suppose for contradiction that $|a|<|s|/2$ and $|b|<|s|/2$. Let $\tilde{s}$ be the largest floating-point number smaller than $|s|$. Then $|a|\leq \tilde{s}/2<s/2$, $|b|\leq \tilde{s}/2 < s/2$, and $|a+b|\leq \tilde{s} < |s|$, a contradiction.
\end{proof}
\begin{lemma}
	Consider Figure \ref{fig:baileyadd}. Let $z_i= (|x_i| \geq |y_i|)\, ?\, x_i : y_i$, $w_i=(|x_i| \geq |y_i|)\, ?\, y_i : x_i$, and $\delta_i=\mathrm{ulp}(z_i)$. Then $\mathrm{uls}(e_{k,j})\geq 2^{-1}\delta_j$.
\end{lemma}
\begin{proof}
	It suffices to consider the smallest $k<j$ with $\mathrm{ulp}(w_k)<\delta_j/2$, i.e., $\mathrm{ulp}(w_k)\leq \delta_j/4$. In this case $\mathrm{ulp}(w_i)\geq \delta_j/2$ for $i <k$, so it suffices to consider $w_i$ with $k\leq i < j$.
	For $k<i\leq j$, $\mathrm{ulp}(w_i) \leq 2^{1-p-2}\delta_j = 2^{-p-1}\delta_j$, i.e., $|w_i|\leq (2^p-1)2^{-p-1}\delta_j<2^{-1}\delta_j\leq \frac 1 2 \mathrm{ulp}(|z_j|)$; for $k\leq i < j$ and $k\leq s \leq j$, $\delta_i \geq 2^{p-1}\delta_j$ and $|w_s|\leq (2^p-1)\delta_j/4 < 2^{p-2}\delta_j \leq 2^{-1}\delta_i$, i.e., $|w_s|<\frac 1 2 \mathrm{ulp}(z_i)$. In other words, except for $f_{k,j-1}+e_{k+1,j}$ in the lower-right corner, all additions in the triangular region above row $k$ (inclusive) and to the left of column $j$ (inclusive) are pass-through:
	$f_{i,s}=w_i$, $e_{i,s}=z_s$ for $( k<i\text{ and } i\leq s \leq j )\text{ or } (k=i \text{ and } i\leq s < j)$. 
	Since $e_{k+1,j}=z_j$, we obtain $\mathrm{uls}(e_{k,j})\geq \delta_j/2$.
	
\end{proof}

This estimate cannot be used directly to prove $\mathrm{uls}(e_{0,j})\geq \mathrm{ulp}(\mathrm{RORS}(e_{0,j+1},e_{0,j+2},\cdots))$, but it can be used to bound some inequalities in the subsequent proof.
\begin{remark}
	The lower bound of $\mathrm{uls}(e_{0,i})$ depends on $\delta_i$, while the upper bound of $\mathrm{ulp}(e_{0,i+1})$ depends on $\delta_j$ with $j\leq i+1$. Hence the two bounds are not enough to determine the relationship between $\mathrm{uls}(e_{0,i})$ and $\mathrm{ulp}(e_{0,i+1})$.
\end{remark} 

\begin{proof}[Proof of theorem \ref{thm:baileyulsulp}]
	
	Since $2^p>2^{n-1}(n+1)+6$, we have $n\leq p-1$.
	
	Let $k = \arg\max_{i} |f_{i,j}|$ and $T_j=|f_{k,j}|$. From $\delta_j/2 \leq \mathrm{uls}(e_{0,j})<\mathrm{ulp}(\mathrm{RORS}(e_{0,j+1},e_{0,j+2},\cdots))$ we obtain an estimate of $T_j$ and of its relation to the sum $\mathrm{RORS}(e_{0,j+1},e_{0,j+2},\cdots)$. 
	
	For $j=0$ it is obvious that $\mathrm{uls}(e_{0,0}) \geq \mathrm{ulp}(\mathrm{RORS}(e_{0,j+1},\cdots e_{0,n}))$, so we assume $1\leq j \leq n-1$.
	
	For $j<t$, $e_{0,t}$ is affected by
	$f_{i,j}$ with $0\leq i \leq j$ and $|f_{i,j}|\leq T_j$, as well as by $x_{i},y_{i}$ with $j<i\leq t$. We have $z_i\leq 2^{(p-1)(i-j-1)}\delta_j$ for $j<i$.
	The contribution of $f_{i,j}$ to $|e_{0,t}|$, which is also that to $e_{i,j+1}$, is at most $T_j(1+\mathbf{u})$, with horizontal distance $t-j-1$ and vertical distance $i$.
	\begin{align*}
	C_1(t) &\leq T_j(1+\mathbf{u}) \Sigma_{i=0}^j C_{t-j-1+i}^{i} \mathbf{u}^{t-j-1}(1+\mathbf{u})^{t-j-1+i} 
	%&\leq T(1+\mathbf{u}) \Sigma_{i=0}^j C_{t-j-1+i}^{i} u^{t-j-1}(1+u)^{t-j-1+j} \\
	%&\leq T(1+\mathbf{u})^t \mathbf{u}^{t-j-1} \Sigma_{i=0}^j C_{t-j-1+i}^{i}\\
	%&=  T(1+\mathbf{u})^t \mathbf{u}^{t-j-1} \frac {(j+1)C_t^{j+1}} {t-j}
	\end{align*}
	\begin{align*}
	\Sigma_{t=j+2}^{n}  C_1(t) &\leq T_j(1+\mathbf{u})\Sigma_{t=j+2}^n \Sigma_{i=0}^j \mathbf{u} C_{t-j-1+i}^i (1+\mathbf{u})^{n-1}\\
	&\leq T_j\mathbf{u}(1+\mathbf{u})^n\left(C_{n+1}^{j+1}-(j+2)\right)\\
	&\leq T_j\mathbf{u}(1+\mathbf{u})^nC_{n+1}^{j+1}
	\end{align*}
	The contribution of $x_i,y_i$ with $j<i$ to $|e_{0,t}|$ corresponds to $e_{i,i}$, where $|e_{i,i}|\leq  2\cdot 2^{-(p-1)(i-j-1)}\delta_j =2^{1-(p-1)(i-j-1)}\delta_j$, with horizontal distance $t-i$ and vertical distance $i$. 
	\begin{align*}
	C_2(t) &\leq \Sigma_{i=j+1}^t 2^{1-(p-1)(i-j-1)}\delta_j C_{t}^{i} \mathbf{u}^{t-i}(1+\mathbf{u})^{t}
	%&=  \delta_j (1+\mathbf{u})^{t} \Sigma_{i=j+1}^t 2^{1-(p-1)(i-j-1)}  C_{t}^{i} \mathbf{u}^{t-i}
	\end{align*}
	\begin{align*}
	\Sigma_{t=j+2}^{n}  C_2(t) &\leq \delta_j  2^{1-(p-1)}\Sigma_{t=j+2}^n (1+\mathbf{u})^t \Sigma_{i=0}^t C_t^i\mathbf{u}^{t-i}\\
	&= 4\delta_j \frac {(1+\mathbf{u})^{2(n+1)}-(1+\mathbf{u})^{2(j+2)}} {2+\mathbf{u}}\\
	\end{align*}
	
	\begin{align*}
	|\mathrm{RORS}(e_{0,j+1},\cdots e_{0,n})| &\leq (1+\mathbf{u})^{n-j-1}\left(T_j(1+\mathbf{u})\Sigma_{i=0}^j(1+\mathbf{u})^{i+1} +2\delta_j(1+\mathbf{u})^{j+1}  + \Sigma_{t=j+2}^n (C_1(t)+C_2(t))\right)\\
	&\leq  (1+\mathbf{u})^{n-j-1}\left(  T_j\left( (1+\mathbf{u})\frac {(1+\mathbf{u})^{j+2}-(1+\mathbf{u})} {\mathbf{u}} +\mathbf{u}(1+\mathbf{u})^nC_{n+1}^{j+1} \right)\right.\\
	&\left.\quad +\delta_j \left( 2(1+\mathbf{u})^{j+1}+4 \frac {(1+\mathbf{u})^{2(n+1)}-(1+\mathbf{u})^{2(j+2)}} {2+\mathbf{u}} \right)\right) \\
	& = T_j\left( \frac {(1+\mathbf{u})^{n+2}-(1+\mathbf{u})^{n-j+1}} {\mathbf{u}} + \mathbf{u}(1+\mathbf{u})^{2n-j-1}C_{n+1}^{j+1} \right) \\
	& \quad + \delta_j \left( 2(1+\mathbf{u})^n + 4\frac {(1+\mathbf{u})^{3n-j+1} -(1+\mathbf{u})^{n+j+3}} {2+\mathbf{u}} \right) \\
	&\leq T_j \left( (j+1) \frac { (1+\mathbf{u})^{n-j+1}} {1-(j+1)\mathbf{u}} +  \mathbf{u}(1+\mathbf{u})^{2n-j-1}C_{n+1}^{j+1} \right) \\
	&\quad + \delta_j \left(2(1+\mathbf{u})^n + 4\mathbf{u}\frac {(2n-2j-2)} {(2+\mathbf{u})(1-(2n-2j-2)\mathbf{u})} \right) \\
	& \leq  T_j \left( j+1 + (j+1)\left(\frac { (1+\mathbf{u})^{n-j+1}} {1-(j+1)\mathbf{u}}-1 \right) + \mathbf{u}(1+\mathbf{u})^{2n-j-1}C_{n+1}^{j+1}   \right) \\
	& \quad + \delta_j \left( 2 + \frac {2n\mathbf{u}} {1-n\mathbf{u}} + \frac {4n\mathbf{u}} {1-2n\mathbf{u}}  \right)\\
	& \leq  T_j \left(  j+1 + (j+1)\left(\frac {\frac{(n-j+1)\mathbf{u}} {1-(n-j+1)\mathbf{u}}+(j+1)\mathbf{u} } {1-(j+1)\mathbf{u}} \right) + \mathbf{u}(1+\mathbf{u})^{2n-j-1}C_{n+1}^{j+1}   \right) \\
	& \quad + \delta_j \left( 2 + 3\gamma_{2n}  \right)\\
	& =  T_j \left(  j+1 + \mathbf{u}(j+1)\left(\frac {j+1+(1-\mathbf{u}(j+1))(n-j+1)} {(1-(j+1)\mathbf{u})(1-(n-j+1)\mathbf{u})}\right) + \mathbf{u}(1+\mathbf{u})^{2n-j-1}C_{n+1}^{j+1}   \right) \\
	& \quad + \delta_j \left( 2 + 3\gamma_{2n}  \right)\\
	& \leq  T_j \left(  j+1 + (j+1)\left(\frac {(n+2)\mathbf{u}} {1-(n+2)\mathbf{u}}\right) + \mathbf{u}(1+\mathbf{u})^{2n-j-1}C_{n+1}^{j+1}   \right) \\
	& \quad + \delta_j \left( 2 + 3\gamma_{2n}  \right)\\
	&\leq T_j\left(j+1+n\gamma_{n+2} +\mathbf{u}(1+\gamma_{2n-2})C_{n+1}^{ j+1})\right) +\delta_j \left( 2 + 3\gamma_{2n}  \right) \\
	& < T_j\left(j+2\right) +\delta_j \left( 2 + 3\gamma_{2n}  \right) 
	\end{align*}
	
	By the inequality in \cite{ErrorCorrect}, Chapter 10, Lemma 7, p.~309,
	\begin{equation*}
	C_n^k \leq \sqrt{\frac {n} {2\pi k(n-k)}} 2^{nH(k/n)},
	\end{equation*}	
	where $H(p)=-p\log_2(p)-(1-p)\log_2(1-p)\leq 1$. Clearly, when $\mod(n,2)=0$, $C_n^k \leq \sqrt{ \frac {n} {2\pi (n/2)^2} } 2^n=\sqrt{\frac {2} {\pi n}}2^n$; when $\mod(n,2)=1$,  $C_n^k \leq \sqrt{ \frac {2n} {\pi (n^2-1)} } 2^n$. Since $\sqrt {\frac {2} {4\pi}}=0.5\times 0.79788456\cdots$ and $\sqrt{\frac{2\cdot 5} {\pi(5^2-1)}}=0.5\times 0.72836562\cdots$, we have $C_n^k \leq 0.8\cdot 2^{n-1}$ for $n\geq 4$. Thus
	\begin{align*}
	n\gamma_{n+2} +\mathbf{u}(1+\gamma_{2n-2})C_{n+1}^{ j+1} & \leq \frac {n(n+2)\mathbf{u}} {1-(n+2)\mathbf{u}} + 2^{-p}\cdot 0.8 \cdot 2^{n} (1+ \frac {(2n-2)\mathbf{u}} {1-(2n-2)\mathbf{u}} )\\
	&\leq  \frac {(p-1)(p+1)2^{-p}} {1-(p+1)2^{-p}} +  \frac {0.4} {1-(2p-4)2^{-p}} \quad (n < p)\\
	& \leq 0.8338983 \quad (\text{maximum attained at } p=7)\\
	& < 1.
	\end{align*}
	
	Alternatively, one may simply use $C_{n+1}^{j+1}<2^{n+1}\leq2^p$, in which case the upper bound on the right-hand side is slightly larger.
	
	Also, $\mathrm{uls}(e_{0,j})\geq \delta_j/2$. To violate the VecSum/Fast2Sum condition, we must have $|\mathrm{ulp}(\mathrm{RN}(\mathrm{RORS}(e_{0,j+1},e_{0,j+2},\cdots e_{0,n})))|\geq \delta_j$, which implies $|\mathrm{RORS}(e_{0,j+1},e_{0,j+2},\cdots e_{0,n})|\geq 2^{p-1}\delta_j(1-\mathbf{u})$.
	\begin{align*}
	T_j &\geq \delta_j \frac {2^{p-1}(1-\mathbf{u}) -2 - 3\gamma_{2n} }  {j+2}\\
	&\geq \delta_j \frac {2^{p-1}-3 }  { j+2}\\
	\end{align*}
	
	For $p\geq 8$ and $n<p$,  $3\gamma_{2n}\leq 3 \frac {2(p-1)2^{-p}} {1-2(p-1)2^{-p}}\leq \frac {21} {121} \leq 1/2$.
	
	From the above we obtain $\delta_j \leq \frac { j+2}  {2^{p-1}-3 } T_j$, and
	\begin{align*}
	|\mathrm{RORS}(e_{0,j+1},\cdots e_{0,n})| 
	%	&\leq T_j\left(j+2\right) +\delta_j \left( 2 + 3\gamma_{2n}  \right) \\
	%	&\leq T_j\left(j+2+ \frac { j+2}  {2^{p-1}-3 } ( 2 + 3\gamma_{2n})\right) \\
	%	&\leq T_j(j+3)
	&\leq T_j \left(j+1+\left (n\gamma_{n+2} +\mathbf{u}(1+\gamma_{2n-2})C_{n+1}^{ j+1}+\frac { j+2}  {2^{p-1}-3 } ( 2 + 3\gamma_{2n})\right)\right)\\
	&\leq T_j\left( j+1+ \left( \frac {(p-1)(p+1)2^{-p}} {1-(p+1)2^{-p}} +  \frac {0.4} {1-(2p-4)2^{-p}}+\right.\right.\\
	&\left.\left.\quad \frac { p}  {2^{p-1}-3 } \left( 2 + 3\frac {(2p-2)\mathbf{u}} {1-(2p-2)\mathbf{u}}\right)\right)\right)  \quad (n < p)\\
	&\leq T_j (j+1+0.81384\cdots)\quad (p=8)\\
	&\leq T_j(j+2)
	\end{align*}
	
	Let $S$ be the largest uls smaller than $\mathrm{ulp}( {\mathrm{RORS}(e_{0,j+1},e_{0,j+2},\cdots)})$. Then $S\leq 2^{-1}\mathrm{ulp}((j+2)T_j)\leq 2^{-p}(j+2)T_j$.
	
	We continue with the propagation argument. Only when a reachable path exists from $x_i,y_i$ to $e_{0,j}$ can a uls of $x_i,y_i$ that is at most $S$ affect $e_{0,j}$; equivalently, $\mathrm{uls}(e_{0,j})$ is at least the minimum of the uls values of the $x_i,y_i$ on which it depends. If along some path from $x_s,y_s$ to $e_{0,j}$ all the uls values of $e_{\cdot,\cdot}$ and $f_{\cdot,\cdot}$ exceed $S$, and all $w_i$ below the path are at least this threshold, then $\mathrm{uls}(e_{0,j})$ is also at least the threshold; the small uls values above the path either cannot propagate to $e_{0,j}$ or are merged into some $f_{\cdot,\cdot}$ on the path, and hence do not lower the estimate of $\mathrm{uls}(e_{0,j})$.
	
	\begin{itemize}
		\item $k=0$: then $|e_{0,j}|\geq 2^pT_j \geq (j+2)T_j\geq |\mathrm{RORS}(e_{0,j+1},e_{0,j+2},\cdots)|$. 
		\item $k=1$: then $\mathrm{ulp}(e_{1,j})\geq 2T_j$ and $\mathrm{uls}(e_{0,j})\geq \frac 1 2 \mathrm{ulp}(e_{1,j})\geq T_j \geq  2S$, or $e_{0,j}=0$.
		\item $k = j$: since $|x_j|+|y_j|\leq 2(2^p-1)\delta_j$ and $T_j=|f_{j,j}|\leq \frac 1 2 \mathrm{ulp}(\mathrm{RN}(|x_j|+|y_j|)) \leq2\delta_j$, this contradicts $T_j\geq \delta_j \frac {2^{p-1}-3 }  { j+2}$.
		\item $2\leq k \leq j-1$: first, $|e_{k,j}|=\mathrm{RN}(e_{k+1,j}+f_{k,j-1})\geq 2^p|f_{k,j}|=2^pT_j  $. By Lemma \ref{lem:binsel}, either $|e_{k+1,j}|\geq 2^{-1}\cdot 2^pT_j$ or $|f_{k,j-1}|\geq 2^{-1}\cdot 2^pT_j$. Suppose the latter; then $|e_{k,j-1}|\geq 2^p\cdot 2^{-1}\cdot 2^pT_j$, so either $|e_{k+1,j-1}|\geq 2^{-1}\cdot 2^p \cdot 2^{-1}\cdot 2^p T_j$ or $ |f_{k,j-2}| \geq 2^{-1}\cdot 2^p \cdot 2^{-1}\cdot 2^p T_j$. Continuing this bisection to the left or upward against the arrows of Figure \ref{fig:baileyadd}, we eventually reach some $e_{s,s}$ with $k<s\leq j$, or some $f_{s,s}$ with $k\leq s<j $. Let the path be $v_0\rightarrow v_1 \rightarrow \cdots \rightarrow v_m$, with starting point $v_0=f_{k,j}$ and endpoint $v_{m}\in\{e_{s,s},f_{s,s}\}$. When $v_i$ ($i>0$) is $e_{t,q}$ with $t\leq q\leq j$, we have $v_{i+1}\in\{e_{t+1,q},f_{t,q-1}\}$ and $|v_{i+1}|\geq 2^{-1}|v_i|$; when $v_i$ ($i>0$) is $f_{t,q}$ with $t \leq q<j$, we have $v_{i+1}\in\{e_{t+1,q},f_{t,q-1}\}$ and $|v_{i+1}|\geq 2^{p-1}|v_i|$. Clearly $m-1=(s-k)+(j-s)=j-k\leq n-1-2=n-3$. By induction, when $v_i$ is $e_{t,q}$, $|e_{t,q}| \geq 2^{-(t-k+j-q)}\cdot 2^{p(j-q+1)}T_j\geq 2^{k-j}2^{p(j-q+1)}T_j$, and when $v_i$ is $f_{t,q}$, $|f_{t,q}| \geq 2^{-(t-k+j-q)}\cdot 2^{p(j-q)}T_j\geq 2^{k-j}2^{p(j-q)}T_j$.
		
		For all $v_i$ with $i>0$, $|v_i|\geq 2^{3-n+p}T_j \geq (j+2)T_j \geq  |\mathrm{RORS}(e_{0,j+1},e_{0,j+2},\cdots)|$, which requires $2^{p}\geq 2^{n-3}(n+1)$. Hence no floating-point number crossing the path from the upper-right side to the lower side can introduce a uls smaller than $S$.
		
		\begin{itemize}
			\item If $v_m = f_{s,s}$ with $s<j$: 
			1. For $i\leq s$: 
			1a. If $s<j-1$, then $|f_{s,s}|\geq 2^{3-n+2p}T_j$, and $\min(\mathrm{uls}(w_i),\mathrm{uls}(z_i))\geq \mathrm{ulp}(w_i) \geq \mathrm{ulp}(w_s) \geq \mathrm{ulp}(2^{3-n+2p}T_j)=2^{3-n+p}T_j\geq 2^3T_j\geq  2^{1-p}(j+2)T_j\geq 2S$.
			1b. If $s=j-1$, then $f_{s,s}$ propagates to $e_{0,j}$ through $e_{k,j}$, and $\mathrm{uls}(e_{k,j})\geq \mathrm{ulp}(e_{k,j})\geq 2|f_{k,j}|=2T_j\geq 2^{1-p}(j+2)T_j=2S$; the case $w_i$ with $i <s$ is as above.
			2. The $w_i$ with $i > s$ lie above the path and cannot affect $e_{0,j}$. 
				\item If $v_m=e_{s,s}$: 
			\begin{itemize}
				\item $s=j$: then $|e_{j,j}|\geq 2^{3-n+p}T_j\geq\frac {2^{3-n+p}(2^{p-1}-3) }  { j+2}\delta_j$, whereas $|e_{j,j}|\leq 2\cdot 2^p\delta_j =2^{p+1}\delta_j$, a contradiction. This requires $2^{3-n+p}(2^{p-1}-3)>(n+1)2^{p+1}$, i.e., $2^{2-n}(2^{p-1}-3)>(n+1)$, i.e., $2^{p}>2^{n-1}(n+1)+ 6$.
				\item $k\leq s < j$: then $z_i\geq 2^{p-1}z_s \geq 2^{p-2}|e_{s,s}|\geq 2^{p-2+3-n+2p}T_j=2^{3p+1-n}T_j$, and $\delta_i \geq 2^{2p+1-n}T_j \geq 2^{1-p}(j+2)T_j \geq 2S$ for $0\leq i < s$. Hence the uls values of the $z_i$ below the path are all at least $2S$, while those above the path have no effect; only the $w_i$ need to be considered.
				
				Suppose the path contains $f_{t,s}$ with $k\leq t < s$. If the first ulp not exceeding $S$ appears in column $q$, i.e., $\mathrm{ulp}(w_i)\geq 2S$ for $i<q$ and $\mathrm{ulp}(w_q)\leq S$, then for $q\leq h,m \leq s$, $w_h\leq w_q \leq 2^pS$, $z_m \geq z_s \geq 2^{-1}|e_{s,s}| \geq 2^{2-n+2p}T_j$, so $w_h/z_m< 2^{n-2-p}S/T_j\leq 2^{n-2-2p}(j+2)\leq 2^{-p-1}$ (this requires $ 2^p\geq 2^{n-1}(n+1)$), i.e., $w_h< \frac 1 2 \mathrm{ulp}(z_m)$. In other words, columns $q$ through $s$ are pass-through, with $f_{h,m}=w_h$ and $e_{h,m}=z_m$. This situation is illustrated in Figure \ref{fig:vmess}.
				\begin{itemize}
					\item $0\leq q \leq t$: then $|f_{t,s}|=w_t \leq w_q \leq 2^pS\leq 2^p(2^{-p}(j+2))T_j=(j+2)T_j$, while the path estimate above gives $|f_{t,s}|\geq 2^{3-n+p}T_j$, a contradiction. This requires $n+1<2^{3-n+p}$, i.e., $2^p>2^{n-3}(n+1)$.
					\item $t < q \leq s$: by the pass-through property above, the relevant values propagate to the upper side of the path and do not pass through $e_{0,j}$.
					\item $s < q$: above the path, so there is no effect on $e_{0,j}$.
					
				\end{itemize}
				\begin{figure}[H]
					\centering
					\includegraphics[width=0.7\linewidth]{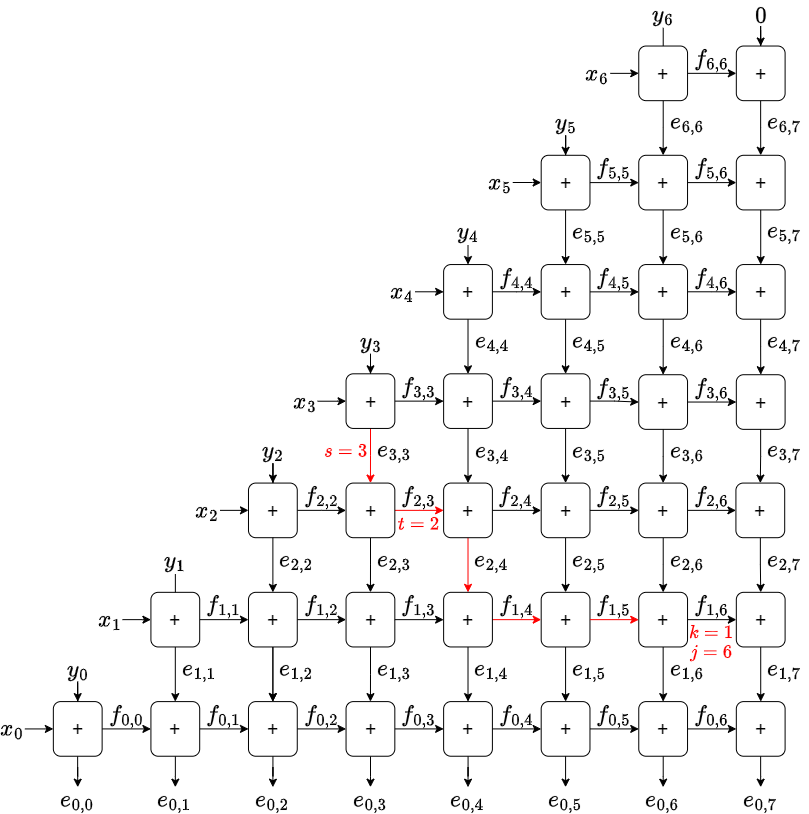}
					\caption{The pass-through region in the case $v_m=e_{s,s}$.}
					\label{fig:vmess}
				\end{figure}

			\end{itemize}
		\end{itemize}
		
	\end{itemize} 
	The additional conditions are, respectively,
	\begin{align*}
	2^p&\geq2^{n-3}(n+1),\\
	2^p&>2^{n-1}(n+1)+6,\\
	2^p&\geq2^{n-1}(n+1),\\
	2^p&>2^{n-3}(n+1).
	\end{align*}
	In fact, only $2^p>2^{n-1}(n+1)+6$ is needed.
	
\end{proof}

\section{Some Remarks}

\subsection{VecSumErrBranch}
It is commonly thought that VecSumErrBranch cannot be vectorized with SIMD because of its conditional branching. In practice, however, this can be achieved with a few programming techniques. While scanning the input sequentially from the beginning, one can choose the next propagated term based on whether the rounding error term returned by Fast2Sum is zero; finally, the redundant interleaving zeros can simply be removed.
\subsection{Merge-based addition}
Similarly, merge-based addition is generally considered difficult to vectorize because the merging step requires conditional branching. However, the VRANGEPD instruction in Intel AVX-512F can directly order two operands by absolute value, so the merging step can be implemented using a simple sorting network with a fixed number of comparisons.

\section{Conclusion}
In this paper, we present a detailed analysis of the behavior of the VecSum algorithm for various input configurations, and then prove the correctness of several renormalization algorithms. We also prove the correctness of the quad-double addition algorithm implemented in the QD library. 

\section{Note on the 2024 Version}
The core work was completed in 2024, and some parts may have been added later. Because the work was not kept under strict version control, I can only confirm that the initial version already contained some key conclusions, such as the correctness of the Bailey addition algorithm. In the years since, I have not had the opportunity to revise it further. Over the past two years, other researchers have reported several new findings, but these have not yet been cited in this note. This version is released as a preprint.
The proofs and methods in this paper were developed without any AI assistance, while the subsequent polishing and editing was assisted by AI tools such as DeepSeek.

	\bibliographystyle{unsrtnat}
	\bibliography{main}  %%% Uncomment this line and comment out the ``thebibliography'' section below to use the external .bib file (using bibtex) .

	
	
	
	
	% \end{thebibliography}

\end{document}